\documentclass[12pt]{article} 
\usepackage{amsmath,amssymb,amsthm,eucal,MnSymbol,amscd,tikz} 
 \usepackage{bbm}  
\usepackage[a4paper,margin=2.5cm,top=1.5in,bottom=3cm]{geometry}
\usepackage{epsf,graphicx}   
\usepackage{hyperref}
\usepackage{float}
\usepackage{tikz}
\usepackage{multirow}
\usepackage{xfrac}
\usepackage{pifont}
\usepackage{color}
\usepackage{tikz}
\usepackage[toc,page]{appendix}
\usetikzlibrary{matrix}
\usepackage{subcaption}
\usepackage[all,cmtip]{xy}
  
\font\tenrm=cmr10

\font\bigss=cmssdc10 scaled 2300

\font\cmsslll=cmss10 at 14 pt  
\newcommand{\e}{\epsilon}

\renewcommand{\gg}{\mathfrak{g}}

\newcommand{\gk}{\mathfrak{k}}  
\newcommand{\gl}{\mathfrak{l}}  
\newcommand{\gm}{\mathfrak{m}}

\newcommand{\gt}{\mathfrak{t}}  
\newcommand{\gu}{\mathfrak{u}}

\newcommand{\so}{\mathfrak{so}}  
\newcommand{\su}{\mathfrak{su}}  
\newcommand{\spin}{\mathfrak{spin}}  
\newcommand{\gsp}{\mathfrak{sp}}

\newcommand\Sp{\mathrm{Sp}}  
\newcommand\GL{\mathrm{GL}}  
\newcommand\SL{\mathrm{SL}}  
\newcommand\SO{\mathrm{SO}}  

\newcommand\SU{\mathrm{SU}}  
\newcommand\U{\mathrm{U}}

\newcommand\Spin{\mathrm{Spin}}  

\DeclareMathOperator\diag{diag}
\DeclareMathOperator\tr{tr}  
\DeclareMathOperator\Tr{Tr}  
\DeclareMathOperator\BD{2D} 
  
\DeclareMathOperator\vol{vol\;}  
  
\DeclareMathOperator\ad{ad}  
\DeclareMathOperator\Hom{Hom}  
  
 \DeclareMathOperator\Ad{Ad}

\DeclareMathOperator\Inn{Inn}

\DeclareMathOperator\Out{Out}
\DeclareMathOperator\Aut{Aut}
\DeclareMathOperator\I{I}
\DeclareMathOperator\Real{Re}

\DeclareMathOperator\Ker{Ker}

\DeclareMathOperator\ord{ord}
\DeclareMathOperator\BT{2T}
\DeclareMathOperator\BO{2O}
\DeclareMathOperator\BI{2I}

\DeclareMathOperator\D{D}
\DeclareMathOperator\T{T}

\DeclareMathOperator\OR{O}

\DeclareMathOperator\Alt{Alt}
\DeclareMathOperator\Sym{Sym}

\DeclareMathOperator\ID{Id}
\DeclareMathOperator\Pin{Pin}
\DeclareMathOperator\Eig{Eig}

\newtheorem{Th}{Theorem}[section]  
\newtheorem{Prop}[Th]{Proposition}  
\newtheorem{Cor}[Th]{Corollary}  
\newtheorem{Lem}[Th]{Lemma}  

\theoremstyle{definition} 
\newtheorem{Def}[Th]{Definition}  
\newtheorem{Ex}[Th]{Example}

\theoremstyle{remark}

\newtheorem*{remark}{Remark}

\newcommand{\bt}{\begin{Th}\ \ }  
\newcommand{\et}{\end{Th}}  
\newcommand{\bp}{\begin{Prop}\ \ }  
\newcommand{\ep}{\end{Prop}}  
\newcommand{\bc}{\begin{Cor}\ \ }  
\newcommand{\ec}{\end{Cor}}  
\newcommand{\bl}{\begin{Lem}\ \ }  
\newcommand{\el}{\end{Lem}}  
\newcommand{\bd}{\begin{Def}\ \ }  
\newcommand{\ed}{\end{Def}}  
\newcommand{\bex}{\begin{Ex}\ \ }  
\newcommand{\eex}{\end{Ex}}  

\newcommand{\pf}{\begin{proof}}
\newcommand{\epf}{\end{proof}}

\newcommand{\be}{\begin{equation}}  
\newcommand{\ee}{\end{equation}}

\newcommand{\arr}{\begin{array}{rlll}}  
\newcommand{\ea}{\end{array}}  
  
\newcommand{\eea}{\end{eqnarray}}  
\newcommand{\bean}{\begin{eqnarray*}}  
\newcommand{\eean}{\end{eqnarray*}}

\begin{document}  

\vskip 1.0 true cm  
\begin{center}  
{\bigss  Isospectral nearly K\"ahler manifolds}  
\vspace{1.5ex}
\vskip 0.5 true cm   
{\cmsslll J.\ J.\ V\'asquez} \\[3pt] 
{\tenrm   Max-Planck-Institut f\"ur Mathematik in den Naturwissenschaften \\ 
Inselstra{\ss}e 22, D-04103 Leipzig, Germany\\
vasquez@mis.mpg.de} 
\end{center}    
\baselineskip=18pt  
\begin{abstract}  
We give a systematic way to construct almost conjugate pairs of finite subgroups of $\Spin(2n+1)$ and $\Pin(n)$ for $n\in \mathbb{N}$ sufficiently large. As a geometric application, we give an infinite family of pairs $M_1^{d_n}$ and $M_2^{d_n}$ of nearly K\"ahler manifolds that are isospectral for the Dirac and Laplace operator with increasing dimensions $d_n>6$. We provide additionally a computation of the volume of (locally) homogeneous six dimensional nearly K\"ahler manifolds and investigate the existence of Sunada pairs in this dimension.
\end{abstract}
\tableofcontents

\section{Introduction}
Nearly K\"ahler manifolds were introduced by A. Gray in connection with the concept of weak holonomy, see \cite{Gr,Gr2}. 
An almost Hermitian manifold $(M,g,J)$ is said to be (strict) \textit{nearly K\"ahler} if $(\nabla_XJ)X=0$ for any vector field $X\in \Gamma(TM)$ and it has no non-trivial local K\"ahler de Rham factor. Here $\nabla$ denotes the Riemannian connection. 

Based on previous results of Clayton and Swann, Nagy showed in \cite{Na} that a complete simply connected nearly K\"ahler manifold is a Riemannian product of twistor spaces of quaternionic K\"ahler manifolds of positive scalar curvature, homogeneous spaces and six-dimensional nearly K\"ahler manifolds. Some other classification results are also known since then. For example, there is a complete list of homogeneous nearly K\"ahler manifolds obtained by Cabrera, Davila and Butruille \cite{B1,B2,GM,GM1}.
 
In spite of the great deal that is known about this class of manifolds, not much information about the spectrum of natural operators defined on them, such as the Laplace or Dirac operator, seems to be available. To the best of the author's knowledge, the majority of results in this vein concern deformations of
geometric structures on nearly K\"ahler manifolds in dimension six: infinitesimal deformations of the nearly K\"ahler structure \cite{MS} and more recently deformations of instantons \cite{CH}. The aim of this work is to give a systematic way to construct pairs of nearly K\"ahler manifolds that are isospectral for the Laplace and Dirac operator. 

Perhaps the most famous criterion for isospectrality is the one of Sunada, see Theorem \ref{sunada}, which states that given a pair of almost conjugate finite subgroups $\Gamma_i$ of the isometry group of a compact manifold $M$, the quotients $M_{\Gamma_i}:=\Gamma_i \backslash M$ are Laplace isospectral on functions. The main result in Section 1 is Theorem \ref{nakayama}, which gives a systematic way to construct almost conjugate subgroups of $\Spin(2n+1)$ and $\Pin(n)$ for $n\in \mathbb{N}$ sufficiently large. The proof of this statement bulids on an idea of M. Larsen given in \cite{L2}. It should be pointed out that computations involving symmetric characters were carried out using MAGMA. Section 2 contains the main result of this paper: 
\bt  \label{isospectral1} There is a strictly increasing sequence of numbers $(d_n)_{_{n\in \mathbb{N}}}$ such that for each $n\in \mathbb{N}$ there is a pair $M^{d_n}_1$ and $M^{d_n}_2$ of non-isometric nearly K\"ahler manifolds that are isospectral for the Dirac and the Hodge-Laplace operator $\Delta^k$ for $k=0,1,...,\dim(M)$.
\et
The scarcity of nearly K\"ahler manifolds in dimension six suggests that the major chance to find isospectral examples in this dimension is by means of locally homogeneous nearly K\"ahler manifolds. The only such manifolds available in the literature are either homogeneous or quotients of $S^3 \times S^3$ with its canonical nearly K\"ahler structure, see  \cite{CV}. After distinguishing the spectrum of the four homogeneous examples in Proposition \ref{VOL} by means of heat invariants and classifying almost conjugate subgroups of $\Spin(4)$ in Proposition \ref{SU2}, it turns out that Sunada's method does not produce such pairs for any choice of almost conjugate subgroups $\Gamma, \Gamma' \subset \Spin(4) \times \{\ID\}$:
\bp If $M_{\Gamma}$ and $M'_{\Gamma'}$ are non-isometric Laplace isospectral locally homogeneous nearly K\"ahler manifolds in dimension six with $\Gamma \subset \I_0^h(M)$ and $\Gamma' \subset \I_0^h(M')$,  then  $M'$ and $M$ are holomorphically isometric to $S^3\times S^3$ endowed with its 3-symmetric nearly K\"ahler structure. Furthermore,
	 Sunada isospectral pairs $M_{\Gamma}$ and $M_{\Gamma'}$ with $M=S^3\times S^3$ and $\Gamma, \Gamma' \subset \SU(2) \times \SU(2) \times \{\ID\}$ are holomorphically isometric.
	\ep
The construction of isospectral pairs of six dimensional locally homogeneous nearly K\"ahler manifolds is a natural step to continue with the investigation of the spectrum of nearly K\"ahler manifolds.\\[0.3cm]
\textit{Acknowledgements:} The author wishes to thank M. Larsen and T. Finis for discussions concerning results in \cite{L1,L2},
as well as A. Adem for pointing out a reference for the group cohomological facts used in the proof of Theorem \ref{nakayama}. He also wants to thank N. Ginoux for several corrections in the previous versions of this paper and G. Weingart for providing help with the algorithmic computations that yield Example \ref{list} and his hospitality during the author's stay in Cuernavaca.\\[0.3cm]
This work was supported by the Max-Planck-institut f\"ur Mathematik in den Naturwissenchaften. 
\section{Group theoretical preliminaries}
\subsection{ Finite subgroups of Spin(4) }\label{subg}
 Any finite subgroup of $\SU(2)$ is conjugate to one of the so-called ADE groups, see e.g. Theorem
1.2.4 in \cite{GAB}. A description of these groups in terms of unit quaternions is given in Table \ref{tabla1}.  
\begin{table}[H]
\begin{center}
  \begin{tabular}{ ||c | c |c | c |c || }
    \hline
		Label & Name &  Order & Generators & Elements \\ \hline
			$\mathbb{A}_{n-1}$ & $\mathbb{Z}_n$ & $n$ & $e^{\frac{i 2 \pi}{n} }$  &$e^{\frac{ 2 \pi i x}{n} } $\\ 
			$\mathbb{D}_{n+2}$ & $\BD_{2n}$ & $4n$ & $j, e^{\frac{i \pi }{n}}$ &$e^{\frac{i  \pi x}{n} }, j e^{\frac{i  \pi x}{n} } $\\ 
			$\mathbb{E}_6$ & $\BT$ & 24 & $\frac{1}{2}(1+i)(1+j), \frac{1}{2}(1+j)(1+i)$& $\BD_{4}\cup \left\{ \frac{\pm 1 \pm i \pm j \pm k}{2}\right\}$ \\ 
			$\mathbb{E}_7$ & $\BO$ & 48 & $\frac{1}{2}(1+i)(1+j), \frac{1}{\sqrt{2}}(1+i)$& $\BT \cup e^{\frac{i\pi}{4}}\BT$\\ 
			$\mathbb{E}_8$ & $\BI$ & 120 & $\frac{1}{2}(1+i)(1+j), \frac{1}{2}(\phi+\phi^{-1}i +j)$& $q^j \BT,  \ j\leq 4$\\ \hline
\end{tabular}
\end{center}
\caption{Finite subgroups of $\SU(2)$.}
\label{tabla1}
\end{table}
Here $n \geq 2$, $\phi= \frac{1+\sqrt{5}}{2}$ denotes the golden ratio and  $q= \frac{1}{2}\left(\phi +\phi^{-1}i +j \right)$.  The subgroups of $\Spin(4)=\SU(2) \times \SU(2)$ are in turn described by means of Goursat's lemma, 
see Theorem 2.1 in \cite{CV}. 
\bl \label{GOURSAT}
Let $G_1, G_2$ be groups. There is a one-to-one correspondence between subgroups $C  \subset G_1 \times G_2$ and quintuples $\mathcal{Q}(C)=\{A,A_0,B,B_0,\theta\}$, 
where $A_0 \triangleleft A \subset G_1$, $B_0 \triangleleft B  \subset G_2 $ and $\theta: \sfrac{A}{A_0} \longrightarrow \sfrac{B}{B_0} $ is an isomorphism.  
\el
A quintuple $(A,A_0,B,B_0,\theta)$ as described in Lemma \ref{GOURSAT} defines a subgroup of $\Spin(4)$ by setting
\begin{align}
 \mathcal{G}(A,A_0,B,B_0,\theta)& = \{ (a,b)\in A\times B| \; \alpha (a) = \beta (b)\} , \label{fiber-prod.1} 
\end{align}
where  
\begin{center}
\begin{tikzpicture}[node distance=2cm, auto]
  \node (A) {$\alpha:A$};
	\node (B) [right of=A] {$\sfrac{A}{A_0}$};
  \node (C) [right of=B]{$\sfrac{B}{B_0}$};  
  \node (D) [right of=C]{ and};
   \node (E) [right of=D]{$\beta:B $};
    \node (F) [right of=E]{$ \sfrac{B}{B_0}$};
  \draw[->] (A) to node {} (B);
  \draw[->] (E) to node {} (F);
	\draw[->] (B) to node {$\theta$} (C);
	  \end{tikzpicture}
\end{center}
are the natural homomorphisms. Conversely, a group $C \subset \Spin(4)$ defines a quintuple $\mathcal{Q}(C)=\{A,A_0,B,B_0,\theta\}$ by setting 
\begin{align}
A= \pi_1 (C) \subset \SU(2) \quad, \quad  B=\pi_2(C)\subset \SU(2), \label{g1} \\
A_0= \pi_1(\Ker(\pi_2|_{_{C}})) , \quad B_0=\pi_2(\Ker(\pi_1|_{_{C}})), \label{g2} 
\end{align}
and $\theta(aA_0)=bB_0$, where $(a,b)\in C$ and $\pi_i:\SU(2) \times \SU(2) \longrightarrow \SU(2) $ denote the projections to the factors.\\[0.3cm]
In order for use to use Lemma \ref{GOURSAT} effectively, we need to describe normal subgroups of ADE groups and the automorphisms of the quotients. The former groups together with the isomorphism type of the corresponding quotients are described in Tables \ref{tabla-p} and \ref{tabla-p.1}, which are borrowed from \cite{FIG}. 

\begin{minipage}{0.5\textwidth}
\begin{flushleft} 
\begin{table}[H]
\begin{center}
  \begin{tabular}{ ||c | c || }
    \hline
		$A_0 \triangleleft A$ 			 & $A/A_0$	 \\ \hline
			$\mathbb{Z}_k \triangleleft \mathbb{Z}_{kl} $ & $ \mathbb{Z}_l $     \\ 
			$\mathbb{Z}_{2k} \triangleleft \BD_{2kl}$ & $\D_{2l}$ 	 \\  
			$\mathbb{Z}_{2k+1} \triangleleft \BD_{2l(2k+1)}$		 & $\BD_{2l}$  			 	 \\ 
			$\mathbb{Z}_{2k+1} \triangleleft \BD_{2(2k+1)}$		 & $\mathbb{Z}_4$  \\
			$\BD_{2k} \triangleleft \BD_{4k}$		 & $\mathbb{Z}_2$  			  \\ 
			$\mathbb{Z}_2 \triangleleft \BT$  & $\T$  \\ \hline
\end{tabular}
\end{center}
\caption{Subgroups  I. }
\label{tabla-p}
\end{table}
\end{flushleft}
\end{minipage}
\begin{minipage}{0.5\textwidth}
\begin{flushright} 
\begin{table}[H]
\begin{center}
  \begin{tabular}{ ||c | c || }
    \hline
		$A_0 \triangleleft A$ 			 & $A/A_0$ 	 \\ \hline
		
$\BD_4 \triangleleft \BT$	 & $\mathbb{Z}_3$      \\ 
		$\mathbb{Z}_2 \triangleleft \BO$ 	 &  $\OR$ 	 \\ 
			$\BD_4 \triangleleft \BO$ 		 &  $\D_6$			 	 \\ 
			$\BT\triangleleft \BO $ 		 &  $\mathbb{Z}_2$			  \\ 
			$\mathbb{Z}_2 \triangleleft \BI$  & $\I$\\ \hline
\end{tabular}
\end{center}
\caption{Subgroups II.}
\label{tabla-p.1}
\end{table}
\end{flushright}
\end{minipage}
\ \\[0.3cm]
Here $\T,\OR, \I$ denote the usual polyhedral groups and $\D_{2n}$ the dihedral group.
\subsection{Automorphisms of quotient groups of ADE groups}\label{automorphisms} 
The present section comprises descriptions of automorphisms groups of quotients of ADE groups that are relevant in the forthcoming sections. The main reference is \cite{CV}, but the material will be however adapted to our needs.

(1) The group of outer automorphisms of $\mathbb{Z}_n$ is given by 
\begin{align}
\Out(\mathbb{Z}_n)&=\{ \varphi(r) \ : \ \gcd (r,n)=1 \} , \nonumber
\end{align}
where $\varphi(r)$ denotes the map $\mathbb{Z}_n \ni x \mapsto x^r \in \mathbb{Z}_n$ in multiplicative notation. 

(2) To describe the outer automorphism group of a dihedral group $\D_{2n}$, consider the following presentation of $\D_{2n}$  
\begin{align}
\D_{2n}&=\left< x,y \ : \ x^2=y^n=(xy)^2=1\right> = \{ y^p  : \ 0 \leq p < n\} \cup \{ xy^p  : \ 0 \leq p < n\} . \nonumber
\end{align}
 Observe $\D_2 = \mathbb{Z}_2$, so we can assume that $n>1$. The case $n=2$ is also special as $\D_4$ is isomorphic to the Klein  Vierergruppe. The automorphism group of $\D_4$ isomorphic to $\Sym(3)$ and acts by permutations of the 3 non trivial involutions. The outer automorphism group of $\D_{2n}$ for $n >2$  is 
		\begin{align}
		\Out (\D_{2n})&=\left< \tau_{a,b} \ : \ (a,b)\in \mathbb{Z}^{\times}_{n} \times \mathbb{Z}_{n} \right> \cong \mathbb{Z}^{\times}_{n} \ltimes \mathbb{Z}_{n},  \nonumber
		\end{align}
	where the action of the affine group $\mathbb{Z}^{\times}_{n} \ltimes \mathbb{Z}_{n}$ on $\D_{2n}$ is given by
	\begin{align}
	\tau_{a,b}(y^p)&=y^{ap} \ \ , \ \ \tau_{a,b}(xy^p)= xy^{ap+b}. \nonumber  
	\end{align}
(Here $\mathbb{Z}_{n}$ denotes the additive group and $\mathbb{Z}_{n}^{\times}$ the multiplicative group of units in the ring $\mathbb{Z}_{n}$.) 
	
(3) Since $\BD_2 \cong \mathbb{Z}_4$,  we consider $\BD_{2n}$ only for $n>1$. 
We have the following  presentation: 
\begin{align}
\BD_{2n}&=\left< s,t  : \ s^2=t^n=(st)^2 \right> =\{ t^p   : \ 0 \leq p <2n \} \cup \{ st^p  \ : \ 0 \leq p <2n \} . \nonumber
\end{align} 
In fact, we can take $s=je^{i\frac{\pi}{n}}$ and $t=e^{i\frac{\pi}{n}}$ when $\BD_{2n}$ is realized as a subgroup
of $\SU (2)$. 
The outer automorphism group of $\BD_{2n}$ for $n>2$ is also an affine  group:
\begin{align}
		\Out (\BD_{2n})&=\left< \tau_{a,b} \ : \ (a,b)\in \mathbb{Z}^{\times}_{2n} \times \mathbb{Z}_{2n} \right> 
		 \cong \mathbb{Z}^{\times}_{2n} \ltimes \mathbb{Z}_{2n}, \nonumber  
		\end{align}
	where the action on $\BD_{2n}$ is given by
	\begin{align}
	\tau_{a,b}(t^p)&= t^{ap} \ \ , \ \	\tau_{a,b}(st^p)= st^{ap+b}. \nonumber
	\end{align}	
	We need to make a distinction for $n=2$. Any automorphism of $\BD_4=\{\pm 1, \pm i , \pm j ,  \pm k\} \subset \SU(2)$ is obtained via conjugation with an element in $\BO$ modulo $\mathbb{Z}_2=\{\pm1\}$. The point-wise action of $\sfrac{\BO}{\mathbb{Z}_2}$ on $\BD_4$ is described below.

\begin{table}[H]
\begin{center}
 \begin{tabular}{ c || c | c | c | c |c || c |c |c}   
						&  $i$ & $j$ &  $k$  & & &  $i$ & $j$ &  $k$ \\ \hline \hline	
			$[i]$		&   $i$ &  $-j$ & $-k$ &      &  $[\frac{1}{\sqrt{2}}(1-i)]$	 & $i $ & $-k $&  $j $ 	\\
		    $[j]$		 & $-i$ & $j$ & $-k$ &        & $[\frac{1}{\sqrt{2}}(j+k)]$	   & $-i $  & $k $&  $j $ 	\\
			$[k]$		 & $-i$  & $-j$ & $k$ 	&     	& $[\frac{1}{\sqrt{2}}(j-k)]$	 & $-i $ & $-k $ & $-j $\\
			 $[\frac{1}{2}(1+i+j+k) ]$  & $ j$  & $ k$ & $ i$ & & $[\frac{1}{\sqrt{2}}(i+k)]$	   & $ k$  & $-j $&  $i $ 	\\
			$[\frac{1}{2}(1-i-j-k)]$   & $ k$  & $ i$& $ j$ & &      $[\frac{1}{\sqrt{2}}(1-k)]$	  & $ -j$  & $i $&  $k $ 	\\

		    $[\frac{1}{2}(1+i-j-k)]$   & $ -j$ & $ k$&  $ -i$	& \quad \quad & 		$[\frac{1}{\sqrt{2}}(i-k)]$	   & $-k $  & $-j $&  $-i $ 	\\
		    $[\frac{1}{2}(1+i+j-k)]$ & $ -k$  & $ i$&  $ -j$ 	& & 		$[\frac{1}{\sqrt{2}}(i+j)]$	  & $j $  & $i $&  $-k $ 	\\	
		    $[\frac{1}{2}(1-i+j-k)]$ & $ -j$  & $ -k$& $ i$ 	& & 		$[\frac{1}{\sqrt{2}}(1+j)]$	  & $-k $  & $ j$&  $i $ 	\\
			$[\frac{1}{2}(1-i-j+k)]$ & $ j$  & $ -k$&  $ -i$  & & 	$[\frac{1}{\sqrt{2}}(1-j)]$    & $k $  & $j $&  $-i $	\\

			$[\frac{1}{2}(1-i+j+k)]$ & $ -k$  & $ -i$& $ j$  & & 	$[\frac{1}{\sqrt{2}}(1+k)]$	  & $ j$  & $-i $&  $k $ \\	

				$[\frac{1}{2}(1+i-j+k)]$   & $ k$  & $ -i$&  $ -j$	& & 	$[\frac{1}{\sqrt{2}}(i-j)]$	   & $-j $  & $-i $&  $-k $	 \\
	
				 $[\frac{1}{\sqrt{2}}(1+i)]$	  & $ i$ & $ k$& $ -j$ 	& & 	&  & &				\end{tabular}
			\end{center} 
			\caption{Action of $\sfrac{\BO}{\mathbb{Z}_2}$ on $\BD_4$.}
 \label{tabla3.1}
			\end{table}
(4) The outer automorphism group of $\BT \subset \SU(2)$ is generated by an involution that exchanges the generators $s=\frac{1}{2}(1+i)(1+j), t=\frac{1}{2}(1+j)(1+i)$, which satisfy the relations $s^3=t^3=(st)^3$. This automorphism is given by conjugation with $\frac{1+j}{\sqrt{2}}\in  \BO \subset \SU(2)$.

(5) The outer automorphism group of $\BO \subset \SU(2)$ is generated by an involution $\varphi$ fixing $s$ and sending $t$ to $-t$, where $s= \frac{1}{2}(1+i+j+k) $ and $t= e^{\frac{i\pi}{4}}$ generate $\BO$.

(6) The outer automorphism group of $\BI \subset \SU(2)$ is generated by an involution $\psi$ which fixes $s$ and sends $t$ to $\frac{- \phi^{-1}-\phi i +k}{2}$ , where $s=\frac{1}{2}(1+i+j+k)$ and $t=\frac{\phi + \phi ^{-1}i+ j}{2}$  generate $\BI$. 
The action of the automorphisms $\varphi \in \Out(\BO)$ and $\psi \in \Out(\BI)$ on conjugacy classes $\mathcal{C}(x)$, for $x \in \BO$ or $x\in  \BI$ respectively, is described in the following tables.

\begin{minipage}{0.48\textwidth}
\begin{flushright} 
\begin{table}[H]
\begin{center}
  \begin{tabular}{ ||c | c | c|| }
    \hline
		Representative 			 & Size & Real parts	 \\ \hline
			$1$ & $1$ & $1$\\ 
			$-1$ & $1$ & $-1$			\\ 
			$s$ 		 &  	$8$ & $\frac{1}{2}$		 			 \\ 
			$t$ 		 &  $6$  & $\frac{1}{\sqrt{2}}$\\ 
			$s^2$ & $8$  & $-\frac{1}{2}$\\ 
			$t^2$ & $6$ & $0$			 \\ 
			$t^3$ 		 & $6$ &$- \frac{1}{\sqrt{2}}$ 			 			 \\ 
			$st$ 		 &  	$12$	 & $0$ \\ \hline
\end{tabular}
\end{center}
\caption{Conjugacy classes in  $\BO$ .}
\label{tabla7}
\end{table}
\end{flushright}
\end{minipage}
\begin{minipage}{0.45\textwidth}
\begin{flushleft} 
\begin{table}[H]
\begin{center}
  \begin{tabular}{ ||c | c | c|| }
    \hline		  
		$\mathcal{C}(x)$ 			 & $\varphi(\mathcal{C}(x))$ & $\Real(\varphi(x))$ 	 \\ \hline
			$\mathcal{C}(1)$ & $\mathcal{C}(1)$ & 1 \\ 
			$\mathcal{C}(-1)$ & $\mathcal{C}(-1)$ & 	-1		\\ 
			$\mathcal{C}(s)$ 		 &  	$\mathcal{C}(s)$ & $\frac{1}{2}$ 		 			 \\ 
			$\mathcal{C}(t)$ 		 & $\mathcal{C}(t^3)$   & $- \frac{1}{\sqrt{2}}$ \\ 
			$\mathcal{C}(s^2)$ & $\mathcal{C}(s^2)$  & $- \frac{1}{2}$ \\
			$\mathcal{C}(t^2)$ & $\mathcal{C}(t^2)$ & $0$ \\
			$\mathcal{C}(t^3)$ & $\mathcal{C}(t)$ & $\frac{1}{\sqrt{2}}$ 			 \\ 
			$\mathcal{C}(st)$ 		 & $\mathcal{C}(st)$ 		 & $0$ \\ \hline
\end{tabular}
\end{center}
\caption{Action of $\varphi$.}
\label{conj2O}
\end{table}
\end{flushleft}
\end{minipage}   

\begin{minipage}{0.48\textwidth}
\begin{flushright} 
	\begin{table}[H] 
\begin{center}
  \begin{tabular}{ ||c | c |c || }
    \hline
		Representative			 & Size 	&   Real parts \\ \hline
			$1$ & $1$ &   $1$ \\ 
			$-1$ & $1$ 			& 	$ -1$ \\ 
			$t$ 		 &  	$12$		 			&	  $\frac{1+\sqrt{5}}{4}$ \\ 
			$t^2$ 		 &  $12$			 			&   $-\frac{1-\sqrt{5}}{4}$  \\ 
			$t^3$ & $12$ &   $\frac{1-\sqrt{5}}{4}$ \\ 
			$t^4$ & $12$ 			& 	$-\frac{1+\sqrt{5}}{4} $ \\ 
			$s$ 		 & $20$ 			 			&	  $\frac{1}{2}$ \\ 
			$s^4$ 		 &  $20$			 			&   $-\frac{1}{2}$  \\ 
			$st$ 		 &  	$30$		 			&   $0$  \\ \hline
\end{tabular}
\end{center}
\caption{Conjugacy classes in  $\BI$. }
\label{tabla6}
\end{table}
\end{flushright}
\end{minipage}
\begin{minipage}{0.45\textwidth}
\begin{flushleft} 
\begin{table}[H]
\begin{center}
  \begin{tabular}{ ||c | c | c|| }
    \hline		   
		$\mathcal{C}(x)$ 			 & $\psi(\mathcal{C}(x))$ & $\Real(\psi(x))$ 	 \\ \hline
			$\mathcal{C}(1)$ & $\mathcal{C}(1)$ & 1 \\ 
			$\mathcal{C}(-1)$ & $\mathcal{C}(-1)$ & -1			\\ 
			$\mathcal{C}(t)$ 		 &  	$\mathcal{C}(t^3)$ & $\frac{1-\sqrt{5}}{4}$ 		 			 \\ 
			$\mathcal{C}(t^2)$ 		 & $\mathcal{C}(t^4)$   & $-\frac{1 + \sqrt{5}}{4}$ \\ 
			$\mathcal{C}(t^3)$ & $\mathcal{C}(t)$  & $\frac{1+ \sqrt{5} }{4}$ \\
			$\mathcal{C}(t^4)$ & $\mathcal{C}(t^2)$ & $- \frac{1- \sqrt{5}}{4}$ \\
			$\mathcal{C}(s)$ & $\mathcal{C}(s)$ & $\frac{1}{2}$ 			 \\ 
			$\mathcal{C}(s^4)$ 		 & $\mathcal{C}(s^4)$ & 	$- \frac{1}{2}$		 			 \\ 
			$\mathcal{C}(st)$ 		 & $\mathcal{C}(st)$ 		 & $0$ \\ \hline
\end{tabular}
\end{center}
\caption{Action of $\psi $. }
\label{conj2I}
\end{table}
\end{flushleft}
\end{minipage} 
\ \\[0.3cm]
(7) The tetrahedral group $\T$ is isomorphic to the alternating group $\Alt(4)$, which has automorphism group $\Sym(4)$, acting by conjugation on the normal
subgroup  $\Alt(4)$. This corresponds to the action of the octahedral group $\OR\cong \Sym(4)$ on its normal subgroup $\T$, which is induced by the action of 
$\BO$ on the normal subgroup $\BT$.  In fact, it can be derived from Table \ref{tabla3.1} that the image of $\OR=\sfrac{\BO}{\mathbb{Z}_2}$ in $\Aut(\T)=\Aut(\sfrac{\BT}{\mathbb{Z}_2})$ is isomorphic to $\Sym(4)$. 

(8) Every automorphism of $\OR$ is inner.  

(9) The outer automorphism group $\Out(\I)$ of the icosahedral group $\I$ is isomorphic to $\mathbb{Z}_2$. The non-trivial outer automorphism of $\BI$ induces an automorphism $\varphi\in \Aut(\I)$ by realizing $\I$ as $\BI/Z(\BI)$, where $Z(\BI)$ denotes the center of $\BI \subset \SU(2)$. Observe now that two elements $[x],[y] \in \SU(2) / Z(\SU(2))$ are conjugate to one another if and only if either $\Real(x)=\Real(y)$ or $\Real(x)=-\Real(y)$. It follows from this observation that $\varphi \in \Aut(\I)$ can not be inner, see Table \ref{conj2I}.\subsection{Almost conjugate subgroups}\label{larsen}
Let $(M,g)$ be a compact Riemannian manifold and $\Gamma_i$ be two finite subgroups acting freely and isometrically on $M$. The classical result of Sunada gives a sufficient condition for the quotients $M_{\Gamma_i}=\Gamma_i \backslash M$ to be isospectral with respect to the submersion metric.
\bt\label{sunada} Let $(M,g)$ be a compact Riemannian manifold and let $G$ be a group acting freely and by isometries on $M$. If there is a bijection $\phi:\Gamma_1 \longrightarrow \Gamma_2$ preserving conjugacy classes of elements in $G$, then the quotients $M_{\Gamma_i}$ are Laplace isospectral. 
\et
Examples of isospectral pairs that are constructed as quotients of a compact manifold $M$ using Sunada's method, which include the ones presented in this work, are in fact strongly isospectral \cite{G}, i.e. 
they are isospectral for the Hodge-Laplace operator $\Delta^k$ for any $k=0,...,\dim(M)$.\\[0.3cm]
If the groups $ \Gamma_i \subset G$ in Sunada's theorem are in fact conjugate in $G$, the situation trivializes: if $a\Gamma_1 a^{-1}=\Gamma_2$, then $a:M \longrightarrow M$ induces an isometry between the finite quotients $M_{\Gamma_i}$. A triple $(G,\Gamma_1,\Gamma_2)$ giving rise to an isospectral pair $M_{\Gamma_i}$ by means of Theorem \ref{sunada} will be called \textit{Sunada triple} and the groups $\Gamma_i$ are said to be \textit{almost conjugate} in $G$.
\bd\label{al.con}
 A pair of finite subgroups $\Gamma_i$ of a Lie group $G$ are \textit{almost conjugate} if there is a bijection $\phi:\Gamma_1 \longmapsto \Gamma_2$ that preserves conjugacy classes of elements in $G$, i.e. any given element $\gamma \in \Gamma_1$ is related to $\phi(\gamma)$ by conjugation with an element $x=x(\gamma) \in G$.  
\ed
\begin{remark} Being almost conjugate is a property on the conjugacy class of subgroups of $G$: a pair of subgroups $\Gamma_i \subset G$ are almost conjugate if and only if $g\Gamma_1g^{-1}$ is almost conjugate to $h\Gamma_2h^{-1}$ for
 any fixed $(g,h)\in G \times G$. 
\end{remark}
Depending on the Lie group $G$ Sunada triples $(\Gamma_1,\Gamma_2,G)$ might or might not be so that $a\Gamma_1a^{-1}=\Gamma_2$ for some $a\in G$.
\bp\label{SU2} The following statements hold.
\begin{itemize}
 \item[(a)] Two almost conjugate subgroups of $\SU(2)$ are conjugate.
 \item[(b)] Two almost conjugate subgroups of $\Spin(4)$ are conjugate.  
 \item[(c)] The groups $\OR(6), \Spin(6)$ and $\Pin(6)$ admit almost conjugate subgroups that are not conjugate. 
\end{itemize}
\ep
\pf Identify $G=\SU(2)$ with the sphere of unit quaternions in the usual way and recall that two elements in $G$ are conjugate to each other precisely if their real parts are the same. In particular, 
if we have two almost conjugate subgroups $\Gamma_i \subset G$, then 
 $\Real(\Gamma_1)=\Real(\Gamma_2)$. This together with $|\Gamma_1|=|\Gamma_2|$ implies (a) in view of the standard ADE classification, see Table \ref{tabla1}. Suppose now that $G=\SU(2) \times \SU(2)$ and let $\Gamma_i \subset G$  be almost conjugate 
subgroups defined by quintuples $(A,A_0,B,B_0,\theta)$ and $(A',A'_0,B',B'_0,\theta')$ respectively, see Lemma \ref{GOURSAT}. \\[0.3cm]
\textbf{Claim.} The pair of almost conjugate groups $\Gamma_i \subset G$ are conjugate to groups defined by quintuples of the form $\mathcal{Q}(\Gamma_i)=(A,A_0,B,B_0, \theta_i)$.
\pf[Proof of the claim.]
The almost conjugacy condition between the groups $\Gamma_i\subset G$ implies that $A_0=\{a \in A \ : \ (a,1) \in \Gamma_1\}$ and $A'_0=\{a \in A' \ : \ (a,1) \in \Gamma_2\}$ are such that $\Real(A_0)=\Real(A'_0)$ and $|A_0|=|A'_0|$. It follows that the groups $A_0,A'_0\subset \SU(2)$ are conjugate in $\SU(2)$ and that 
$$|A|=|\pi_1(\Gamma_1)|=\frac{|\Gamma_1|}{|A_0|}=\frac{|\Gamma_2|}{|A'_0|}=|\pi_1(\Gamma_2)|=|A'|. $$  
Almost conjugacy of the groups $\Gamma_i \subset G$ implies that $\Real(A)=\Real(A')$, and hence that $A\subset \SU(2)$ and $A'\subset \SU(2)$ are conjugate in $\SU(2)$. Consequently, we can assume that $\Gamma_1=\mathcal{G}(A,A_0,B,B_0,\theta_1)$ and $\Gamma_2=\mathcal{G}(A,A'_0,B,B'_0,\theta_2)$, see Section \ref{subg} where we introduced the latter notation. \\[0.3cm]
The groups $A_0$ and $A'_0$ are isomorphic normal subgroups of a given ADE group $A \subset \SU(2)$ satisfying the relation $\Real(A_0)=\Real(A'_0)$. 
The following assertion shows that we can assume that $A_0=A'_0$ without loss of generality, and hence the claim.\footnote{See Table \ref{tabla-p} and \ref{tabla-p.1} for information on normal subgroups of ADE groups.} \epf 

\textbf{Assertion.} Two isomorphic normal subgroups $A_0,A'_0 \triangleleft A$ of a given ADE group $A \subset \SU(2)$ that satisfy $\Real(A_0)=\Real(A'_0)$ must be equal, unless $A=\BD_4$ and $A_0 \cong \mathbb{Z}_4$. In the latter case, the groups $A_0$ and $A'_0$ are exchanged by conjugation with an element $g\in \SU(2)$ which leaves $A$ invariant.
\pf[Proof of the assertion] Two isomorphic normal subgroups of either $\BT,\BO$ or $\BI$ are computationally verified to be equal. Let $A=\mathbb{Z}_{mn}=\left<e^{\frac{2\pi i }{mn}}\right>$ and $A_0\triangleleft A$ be isomorphic to $\mathbb{Z}_m$. Each normal subgroup of $A \subset \SU(2)$ isomorphic to $\mathbb{Z}_m$ induces a surjective homomorphism $A \longrightarrow \mathbb{Z}_n$. However, there are exactly $\Phi(n)$ such homomorphisms $\alpha(r)(e^{\frac{2\pi i }{mn}})=e^{\frac{2\pi i r}{n}}$, where $r\in \mathbb{Z}_n^{\times}$ and $\Phi$ is the Euler function. These homomorphisms  have kernel exactly $\mathbb{Z}_m=\left<e^{\frac{2\pi i }{m}}\right>$.\\[0.3cm]
We are left with the case $A=\BD_{2l}$ for $l\geq 2$. The group $A=\BD_4$ has exactly four normal subgroups. Three of them are isomorphic to $\mathbb{Z}_4$ and generated by the imaginary units $ i,j,k$, and the remaining normal subgroup which is its center. The three normal subgroups with four elements are related by conjugation with an element in $\BO$, see Table \ref{tabla3.1}.\\[0.3cm]
Let us suppose that $A=\BD_{2l(2k+1)}$ for $k\geq 1$, $l>1$ and that $A_0\triangleleft A$  is isomorphic to $\mathbb{Z}_{2k+1}$.
We assert that $A_0=\left< e^{\frac{2 \pi i }{2k+1}}\right>$. For this we just need to consider homomorphisms $\BD_{2l(2k+1)} \longrightarrow \BD_{2l} \cong A/A_0$. It is not difficult to see that there are exactly $2l\Phi(2l)$ surjective homomorphisms $\alpha(x,r):A \longrightarrow \BD_{2l}$, where $0 \leq x< 2l$ is an integer, $r\in \mathbb{Z}_{2l}^{\times}$ and $\Phi$ is the Euler function. The automorphisms $\alpha(x,r)$ map the generators of $A=\BD_{2l(2k+1)}$ as follows 
$$\alpha(x,r)(j)= je^{\frac{i\pi x}{l}}, \quad \alpha(x,r)\left(e^{\frac{i \pi }{l(2k+1)}}\right)= e^{\frac{i\pi r}{l}}. $$
The condition $$\alpha(x,r)\left(e^{\frac{i \pi y}{l(2k+1)}}\right)=1$$ implies that $y= 0 \mod l$, and so we get  
$\Ker(\alpha(x,r))=\left< e^{\frac{2 \pi i }{2k+1}}\right>.$
Let us suppose that $l=1$ and that $\alpha:\BD_{2(2k+1)}\longrightarrow \BD_2$ is a surjective homomorphism with $\Ker(\alpha)$ isomorphic to $\mathbb{Z}_{2k+1}$. Since the order $\ord(je^{\frac{i\pi x}{2k+1}})=4$ for any integer $x\in \mathbb{Z}$, the condition $\Ker(\alpha) \cong \mathbb{Z}_{2k+1}$ implies that no element of this form can be mapped to the trivial element. So we must have that $\Ker(\alpha) = \left< e^{\frac{2 \pi i }{2k+1}}\right>$.\\[0.3cm]
Let us consider $A=\BD_{2kl}$ and $A_0 \triangleleft A $ to be isomorphic to $ \mathbb{Z}_{2k}$. Surjective homomorphisms from
$A=\BD_{2kl}$ to $\D_{2l}$ are of the form $\alpha(x,r)$ for an integer $0 \leq x < l$ and a unit $r \in \mathbb{Z}_l^{\times}$. Those automorphisms are given by
$$\alpha(x,r)(j)= je^{\frac{i\pi x}{l}}\mathbb{Z}_{2k}, \quad \alpha(x,r)\left(e^{\frac{i \pi }{lk}}\right)= e^{\frac{i\pi r}{l}}\mathbb{Z}_{2k}, $$
where we realized $\D_{2l}$ as the quotient $\BD_{2kl}/\mathbb{Z}_{2k}$.
Just as before 
$\Ker(\alpha(x,r))=\left< e^{\frac{2 \pi i }{2k}}\right>,$
and so $A_0=\left< e^{\frac{ \pi i }{k}}\right>$. Lastly, let $A=\BD_{4k}$ and suppose that $A_0 \triangleleft A$ is isomorphic to $ \BD_{2k}$. The quotient $A/A_0$ is isomorphic to $\mathbb{Z}_2$, and so we must consider surjective homomorphisms  
$\BD_{4k} \longrightarrow \mathbb{Z}_2$. These are of the form $\alpha(r,s)$ and their action on the generators of $\BD_{4k}$ is as follows
$$\alpha(r,s)(j)= r, \quad \alpha(r,s)\left(e^{\frac{i\pi }{2k}}\right)= s, $$
where $r,s\in \mathbb{Z}_2=\{1,-1\}$.
There is only one of these homomorphisms having kernel isomorphic to $\BD_{2k}$. \epf 
We consider consider without loss of generality that the almost conjugate subgroups $\Gamma_i\subset G$ are defined by quintuples of the form $\mathcal{Q}(\Gamma_i)=(A,A_0,B,B_0,\theta_i)$. In order to show claim (b) we proceed considering the following cases separately 
\begin{itemize}
 \item[(I)] $ A_0=A$ and $B_0=B$,
 \item[(II)] $A_0=B_0=\{1\}$,
 \item[(III)] $A_0=\{1\}$ and $B_0 \neq\{1\}$,
 \item[(IV)] $B_0 \neq\{1\}$ and $A_0 \neq\{1\}$.
\end{itemize}
\textbf{Case I.} The Subgroups $\Gamma_i \subset G$ considered here are of the form $A_i \times B_i$ and clearly any two such subgroups that are almost conjugate must be conjugate.\\[0.3cm]
\textbf{Case II.} $A_0=B_0=\{1\}$.\\[0.3cm]
We proceed by considering all possibilities for $A\subset \SU(2)$.
\begin{enumerate}
\item Let $\underline{A=\mathbb{Z}_n}$. The groups $\Gamma_i \subset G$ are of the form 
$$\Gamma_j= \{ (e^{\frac{2\pi i x}{n}},e^{\frac{2\pi i r_jx}{n}})\ : \ x \in \mathbb{Z}\}, \quad \quad r_j\in \mathbb{Z}_n^{\times}.$$
The almost conjugacy condition between the groups $\Gamma_j \subset G$ reads as follows: for any integer $x\in \mathbb{Z}$ there is $y=y(x) \in \mathbb{Z}$ such that 
$$\cos\left(\frac{2\pi r_1x}{n}\right)=\cos\left(\frac{2\pi r_2y}{n}\right)\quad \quad , \quad \quad \cos\left(\frac{2\pi x}{n}\right)=\cos\left(\frac{2\pi y}{n}\right) .$$
It follows that $r_2=\varepsilon r_1 \mod n$ for some value $\varepsilon \in \{1,-1\}$. Conjugation with e.g. $(1,j)\in G$ exchanges the groups $\Gamma_1$ and $\Gamma_2$.
\item Let $\underline{A=\BD_{2n}}$. That is, we consider almost conjugate subgroups of the form
$$\Gamma_i=\{(e^{\frac{\pi i x}{n}},e^{\frac{\pi i a_kx}{n}}), (je^{\frac{\pi i x}{n}},je^{\frac{\pi i (a_kx+b_k)}{n}})\ : \ a_k \in \mathbb{Z}_{2n}^{\times}, \ b_k \in \mathbb{Z}, \ x \in \mathbb{Z}\}.$$
If the groups $\Gamma_i \subset G$ are almost conjugate, then a similar reasoning as in (1) implies that either $a_2=\varepsilon a_1 \mod 2n$ for some $\varepsilon \in \{1,-1\}$. In case $\varepsilon =1$, then the groups $\Gamma_i \subset G$ are related by conjugation with the element $(1,e^{\frac{i\pi (b_1-b_2)}{2n}}) \in G$, whereas if $\varepsilon =-1$
conjugation with $(1,je^{\frac{i\pi (b_1+b_2)}{2n}})\in G$ will exchange them.
\item Let $\underline{A\in \{ \BT,\BO,\BI\}}$ and recall that the outer automorphism group of $A$ is generated by a single involution $\varphi\in \Out(A)$, see Section \ref{automorphisms}. For $A=\BT$ this involution is obtained
by conjugation with an element in $\SU(2)$, which is not the case for $A=\BO, \BI$. Observe there is an element $a\in A$ such that $\varphi(a)\in \SU(2)$ is not conjugate to $a\in \SU(2)$, see Tables \ref{conj2I} and \ref{tabla7}. 
In particular, if the groups $\Gamma_i= \mathcal{G}(A,(1),A,(1),\theta_i)$ are almost conjugate, then either $\theta_i\in \Inn(\SU(2))$ or $\theta_i \in  \Inn(\SU(2)) \cdot \varphi$. In both cases the resulting groups
are easily checked to be conjugate.
\end{enumerate}
As $\mathbb{Z}_2$ has no non trivial outer automorphisms, any pair of almost conjugate subgroups $\Gamma_i=\mathcal{G}(A,A_0,B,B_0,\theta_i)$ such that the quotient 
$A/A_0$ is isomorphic to $ \mathbb{Z}_2$, must be conjugate. Such cases require don't require further attention, and so they will be omitted in the sequel. \\[0.3cm]
\textbf{Case III.} Suppose that $A \triangleright A_0 \supset  \mathbb{Z}_2$ and $B_0=(1)$.\\[0.3cm]
We analyze case by case for $B \subset \SU(2)$ according Tables \ref{tabla-p} and \ref{tabla-p.1}.
\begin{enumerate}
	\item Let $\underline{B=\mathbb{Z}_k}$. Suppose the groups $\Gamma_i= \mathcal{G}(\mathbb{Z}_{kl},\mathbb{Z}_l,\mathbb{Z}_k,(1),\theta(r_i))$ are almost conjugate for some choice of units $r_i\in \mathbb{Z}_k^{\times}$. A similar argument as in Case II (1) reveals that $r_2=\varepsilon r_1 \mod k$ for some $\varepsilon \in \{ 1,-1\}$. So the groups $\Gamma_i \subset G$ are conjugate to each other.
	\item Let $\underline{B=\BD_{2l}}$ and suppose that $\underline{l>1}$. Similarly as in Case II (2), almost conjugate subgroups of the form 
$$\Gamma_i=\mathcal{G}(\BD_{2l(2k+1)},\mathbb{Z}_{2k+1},\BD_{2l}, (1), \theta_{a_i,b_i})\subset G$$ are such that $a_2=\varepsilon a_1 \mod 2l$ for some 
     $\varepsilon \in \{ 1,-1\}$, and hence conjugate. The case in which $\underline{l=1}$ stands apart, since $\BD_2=\mathbb{Z}_4$. We can choose $\BD_2= \left< j \right> \subset \SU(2)$ and observe that the groups $\Gamma_i$ are conjugate 
     via an element $(e^{\frac{ \pi i }{2(2k+1)}},w)$, where $w \in {\BO}$ is an element such that $c_w(j)=-j$, see Table \ref{tabla3.1}.
		\item The cases involving $\underline{B=\mathbb{Z}_3}$ will be analyzed in item (3) of Case IV.
\end{enumerate}
\textbf{Case IV.} $B_0 \neq\{1\}$ and $A_0 \neq\{1\}$.\\[0.3cm]
Here we discriminate according the isomorphism type of $F =A/A_0$. 
\begin{enumerate}
\item Let $\underline{F \in \{ \T, \OR \}}$. The automorphisms of $F\cong \sfrac{A}{\mathbb{Z}_2}$ are induced by conjugation with an element in $\SU(2)$ for $A\in \{\BO,\BT\}$. In consequence, 
a pair of subgroups of the form  $ \Gamma_i=\mathcal{G}(A, \mathbb{Z}_2,A, \mathbb{Z}_2,\theta_i) $ with $ \theta_i \in \Aut(F)$, are conjugate.
\item Let $\underline{F=\I}$ and let $\Gamma_i=\mathcal{G}(\BI,\mathbb{Z}_2,\BI,\mathbb{Z}_2,\theta_i)$ be a pair of almost conjugate subgroups of $G$. Let us consider the standard representation $\I \subset \SO(3)$ and recall that for the non-trivial outer automorphism $\varphi\in \Out(\I)$ there an element $a \in \I$ which is not conjugate to $\phi(a) \in \I$ in $\SO(3)$, see Table \ref{conj2I} and item (9) of Section \ref{automorphisms}. This observation allows us to conclude that a necessary condition to make the groups $\Gamma_i\subset G$ almost conjugate, is to have either $\theta_i \in \Inn(F)$ or $\theta_i\in   \Inn(F) \cdot \varphi$. In both cases the groups $\Gamma_i$ are readily seen to be conjugate. 
\item Let $\underline{F = \mathbb{Z}_3}$. We consider first the quintuples $(A,A_0,\mathbb{Z}_{3n},\mathbb{Z}_n,\theta)$ for $n\geq 1$. Groups of the form 
$\Gamma_i = \mathcal{G}(\mathbb{Z}_{3n}, \mathbb{Z}_n,\BT,\BD_4,  \theta_i) $, where $\theta_i \in \Aut(\mathbb{Z}_3)$, are given by
\begin{align} 
	\Gamma_1= \biggl\{(e^{\frac{2 \pi i x}{n}},z_0),(e^{\frac{2\pi i(1+3y) }{3n} },z_1), (e^{\frac{2 \pi i(2+3z)  }{3n} },z_2) \ : \  z_0 \in \BD_4, \  z_1 \in  \frac{1-i-j-k}{2} \BD_4, \nonumber 
\\ z_2 \in  \frac{1+i+j+k}{2} \BD_4, \  x,y,z \in \mathbb{Z} \biggl\} \nonumber
     \end{align}
		and
     \begin{align} 
	\Gamma_2= \biggl\{(e^{\frac{2 \pi i x}{n}},z_0),(e^{\frac{2\pi i(1+3y) }{3n} },z_1), (e^{\frac{2 \pi i(2+3z)  }{3n} },z_2) \ : \ z_0 \in \BD_4, \  z_1 \in  \frac{1+i+j+k}{2} \BD_4, \nonumber 
\\ z_2 \in  \frac{-1+i+j+k}{2} \BD_4, \  x,y,z \in \mathbb{Z} \biggl\} .\nonumber
     \end{align}
     The groups $\Gamma_i$ are conjugate to each other via elements $(1,w)$, where $w \in \BO$ is so that the conjugation mapping  $c_w \in \Aut(\SU(2))$ satisfies the following conditions:
     $$c_w \left(  \frac{1+i+j+k}{2} \BD_4\right)= \frac{-1+i+j+k}{2} \BD_4, $$ 
     $$ c_w \left(  \frac{-1+i+j+k}{2} \BD_4\right)= \frac{1+i+j+k}{2} \BD_4,$$
see Table \ref{tabla3.1}.

The groups $\Gamma_i=\mathcal{G}(\mathbb{Z}_{3k},\mathbb{Z}_k,\mathbb{Z}_{3k'},\mathbb{Z}_{k'},\theta_i)$ 
with $\theta_i\in \Aut(\mathbb{Z}_3)$ are not almost conjugate unless $\theta_1=\theta_2 \in \Aut(\mathbb{Z}_3)$, case in which they are indeed conjugate. This is simply because elements of the form 
$$ (e^{\frac{2 \pi i (1+3x)}{3k}},e^{\frac{2 \pi i (1+3y)}{3k'}}), \quad x,y \in \mathbb{Z}$$
can not be conjugate to either 
$$ (e^{\frac{2 \pi i (2+3x_1)}{3k}},e^{\frac{2 \pi i (1+3y_1)}{3k'}}) \quad \text{or} \quad (e^{\frac{2 \pi i (1+3x_1)}{3k}},e^{\frac{2 \pi i (2+3y_1)}{3k'}})$$
for any $x_1,y_1 \in \mathbb{Z}$.  Subgroups of the form $\Gamma_i=\mathcal{G}(\BT,\BD_4,\BT,\BD_4,\theta_i)\subset G$ are checked to be conjugate via an element of the form $(w,w')\in \BO \times \BO \subset G$. 
\item Let $\underline{F = \mathbb{Z}_4}$. Just as in Case III (1), one checks that almost conjugate subgroups $\Gamma_i\subset G$ of the form $\Gamma_i=\mathcal{G}(\mathbb{Z}_{4k},\mathbb{Z}_{k}, \mathbb{Z}_{4k'},\mathbb{Z}_{k'},\theta_i)$ must be conjugate. Let us consider subgroups $\Gamma_i= \mathcal{G}(\mathbb{Z}_{4k},\mathbb{Z}_k,\BD_{2(2k'+1)},\mathbb{Z}_{2k'+1},\theta_i)$ with $\theta_i \in \mathbb{Z}_4^{\times}= \mathbb{Z}_2$. The groups in question are given by 
\begin{align} 
	\Gamma_1= \biggl\{  (e^{\frac{2\pi i x_0}{k}},e^{\frac{2 \pi i y_0}{2k'+1}}),\ (e^{\frac{i \pi (1+4x_1)}{2k}},e^{\frac{i \pi (1+2y_1)}{2k'+1}} ), \ (e^{\frac{i \pi (1+2x_2)}{k}},je^{\frac{2i \pi y_2}{2k'+1}} ),  \nonumber 
\\  (e^{\frac{i \pi (3+4x_3)}{2k}},je^{\frac{i \pi (2y_3+1)}{2k'+1}} ) \ : \ x_j,y_j \in \mathbb{Z}\biggl\} ,\nonumber
     \end{align}
\begin{align} 
	\Gamma_2= \biggl\{  (e^{\frac{2\pi i x_0}{k}},e^{\frac{2 \pi i y_0}{2k'+1}}),\ (e^{\frac{i \pi (1+4x_1)}{2k}}, je^{\frac{2i \pi y_1}{2k'+1}} ), \ (e^{\frac{i \pi (1+2x_2)}{k}}, e^{\frac{i \pi (1+2y_1)}{2k'+1}} ),  \nonumber 
\\  (e^{\frac{i \pi (3+4x_3)}{2k}}, je^{\frac{i \pi (2y_3+1)}{2k'+1}}) \ : \ x_j,y_j \in \mathbb{Z}\biggl\} .\nonumber
     \end{align}
The almost conjugacy condition implies that for any $x\in \mathbb{Z}$ there is a $y \in \mathbb{Z}$ so that 
$$\cos\left(\frac{\pi x}{2k}\right)=\cos\left(\frac{\pi y}{2k}\right).$$
This implies that for some value $\varepsilon \in \{ 1,-1\}$ we have 
\begin{align}
x&=\varepsilon y \mod 4. \label{formulilla}
\end{align}
 In particular, the groups $\Gamma_i\subset G$ are not almost conjugate, for if they were, then the elements 
$$(e^{\frac{i \pi (1+4x)}{2k}},e^{\frac{i \pi (1+2y)}{2k'+1}} )\in \Gamma_1, \quad \ (e^{\frac{i \pi (1+2z)}{k}}, e^{\frac{i \pi (1+2w)}{2k'+1}} )\in \Gamma_2$$
have to be conjugate to one another for some value of $x,y,z,w\in \mathbb{Z}$, which is clearly impossible in view of identity \eqref{formulilla}. Let 
$$\Gamma_i = \mathcal{G}(\BD_{2(2k+1)},\mathbb{Z}_{2k+1},\BD_{2(2k'+1)},\mathbb{Z}_{2k'+1},\theta(r_i))\subset G,$$
where $r_1=1$ and $r_2=3$. These groups are not almost conjugate, as one can check there is an element $(x,y)\in \Gamma_1$ such that $\Real(x)=\Real(y)=0$, but there is no such element in $\Gamma_2$.
\item Let $\underline{F= \mathbb{Z}_l}$ for $\underline{l \neq 2,3,4}$. Consider $\Gamma_i=\mathcal{G}(\mathbb{Z}_{kl},\mathbb{Z}_k,\mathbb{Z}_{pl},\mathbb{Z}_p,\theta(r_i))\subset G$ with $r_i \in \mathbb{Z}_l^{\times}$. In this case 
$$\Gamma_i = \{ (e^{\frac{2\pi i x}{kl}},e^{\frac{2\pi i}{pl}(r_ix+ly)})\ : \ x,y \in \mathbb{Z}\}.$$
The almost conjugacy condition between the groups $\Gamma_i$ reads as follows: for any $x,y\in \mathbb{Z}$, there are integers $x'=x'(x,y)$ and $y'=y'(x,y)$ such that 
$$\cos\left( \frac{2\pi x}{kl}\right) = \cos\left( \frac{2\pi x'}{kl}\right), \ \cos\left(\frac{2\pi }{pl}(r_1x+ly) \right)=\cos\left(\frac{2\pi }{pl}(r_2x'+ly') \right).$$
Equivalently 
\begin{align}
r_1x+ly&= \varepsilon r_2x'+ \varepsilon ly' \mod pl, \label{uno} \\
 x&=\varepsilon' x' \mod kl, \label{dos}
\end{align}
for some $\varepsilon, \varepsilon' \in \{1,-1\}$. From equation \eqref{uno}, we must have that 
\begin{align}
r_1x&=\varepsilon r_2x' \mod l. \label{tres}
\end{align}
Equation \eqref{dos} and \eqref{tres} imply that either 
$r_1=r_2 \mod l$ or $r_1=-r_2 \mod l.$
 In the former case, the groups $\Gamma_i$ are the same, whereas in the latter they are related by conjugation, e.g. with an element $(1,j)\in G$. 
\end{enumerate}
This finishes the proof of claim (b). At last, for claim (c) consider the following almost conjugate diagonal subgroups of $\SO(6)$, which correspond to certain linear codes on $\mathbb{Z}_2^6$ that appear 
first in \cite{Con}
\begin{align} 
	\Gamma_1= \{  (1,1,1,1,1,1), (-1,-1,-1,-1,-1,-1), (-1,-1,1,1,1,1) , \nonumber 
\\   (-1,1,-1,1,1,1), (1,-1,-1,1,1,1), (-1,1,1,-1,-1,-1), \label{grupo1} \\
(1,-1,1,-1,-1,-1), (1,1,-1,-1,-1,-1)\} , \nonumber
     \end{align}
\begin{align} 
	\Gamma_2= \{  (1,1,1,1,1,1), (-1,-1,-1,-1,-1,-1), (-1,-1,1,1,1,1) , \nonumber 
\\   (1,1,-1,-1,1,1), (1,1,1,1,-1,-1), (-1,-1,-1,-1,1,1), \label{grupo2} \\
,(-1,-1,1,1,-1,-1), (1,1,-1,-1,-1,-1)\} ,\nonumber
     \end{align}
     It was noted in Example 2.4 of \cite{RSW}, that these subgroups are not conjugate in $\SO(6)$. The groups given in \eqref{grupo1} and \eqref{grupo2} are also not conjugate in $\OR(6)$. To see this,  write $\OR(6) =\SO(6) \rtimes \mathbb{Z}_2$, 
where the $\mathbb{Z}_2$-action on $\SO(6)$ is given by conjugation with the diagonal matrix $(-1,1,1,1,1,1)$. Since the given 
$\mathbb{Z}_2$-action leaves the groups $\Gamma_i\subset \SO(6)$ invariant, these groups are conjugate in $\OR(6)$ precisely when they are so in $\SO(6)$, which is not the case.
Let $\pi:\Spin(6) \longrightarrow \SO(6)$ be the standard covering and consider the groups  
$\overline{\Gamma}_i=\pi^{-1}(\Gamma_i)\subset \Spin(6)$. These groups $\overline{\Gamma}_i\subset \Spin(6)$ have the same cardinality and can be easily seen to be almost conjugate in $\Spin(6)$. 
If $\overline{ \Gamma}_1$ and $\overline{\Gamma}_2$ were related by conjugation in $\Pin(6)$, then $\Gamma_i \subset \SO(6)$ would be conjugate in $\OR(6)$, which is again not the case.
 \epf
The case in which the bijection $\phi:\Gamma_1 \longrightarrow \Gamma_2$ between the subgroups $\Gamma_i$ of $G$ in Definition \ref{al.con} is indeed an isomorphism was (partly) addressed by M. Larsen for various
Lie groups \cite{L1,L2} and it has connections with the theory of 
automorphic forms \cite{BL}. 
\bd \label{sl} Let $G$ be a Lie group and let $\Gamma$ be a finite group. Two homomorphisms $\phi_i:\Gamma \longrightarrow G$ are \textit{element-wise} or \textit{almost} conjugate,
if for any $\gamma \in \Gamma$ there is an element $x(\gamma) \in G$ such that 
\begin{align}
\phi_1(\gamma)=x(\gamma)\phi_2(\gamma)x(\gamma)^{-1}. \label{achom}
\end{align}
In addition, the quadruple $(G,\phi_1,\phi_2,\Gamma)$ is called \emph{admissible} if the element-wise conjugate homomorphisms $\phi_i:\Gamma \longrightarrow G$ are so that the groups $\phi_i(\Gamma)\subset G$ are not related by any inner automorphism of $G$.
\ed
\begin{Ex} There are admissible quadruples of the form $(\Spin(n),\phi_1,\phi_2,\Sym(10))$ for $n\geq 35$, see Proposition 2.6 in \cite{L2}.
\end{Ex}
We aim to give a procedure to construct pairs of almost conjugate subgroups of $\Spin(n)$ and $\Pin(n)$ for suitable $n >6$, by producing admissible quadruples of the form $(\Spin(n),\phi_1,\phi_2,\Sym(m))$ and $(\Pin(n),\phi_1,\phi_2,\Sym(m))$ in a systematic manner. To prove this statement, Theorem \ref{nakayama}, we need some preliminary steps.
\begin{Lem}\label{lift} Let $\Gamma=\Sym(m)$ and let $\rho:\Gamma \longrightarrow \OR(n)$ be a representation with character $\chi$. Denote by 
$\pi: \Spin(n) \longrightarrow \SO(n)$ the standard two-fold covering and let $x_1,x_2 \in \Gamma$ be any disjoint transpositions. The following statements hold. 
\begin{itemize}
 \item[(a)] $\rho$ takes values in $\SO(n)$ if and only if $ \chi(x) = n \mod 4$.
 \item[(b)] Let $x\in \{x_1,x_1x_2\}$ and suppose that $\rho(x) \in \SO(n)$. The element $\rho(x)$ lifts to an element of order two in $\Spin(n)$ if and only if 
 $\chi(x) = n \mod 8$.
\item[(c)] Let $x\in \{x_1,x_1x_2\}$ and suppose that $\rho(x) \in \SO(n)$. The element $\rho(x)$ lifts to an element of order four in $\Spin(n)$ if and only if 
 $\chi(x) = n \mod 4$ and  $\chi(x) \neq n \mod 8$. 
\end{itemize}
\end{Lem}
\pf Let $x\in \Gamma=\Sym(m)$ be either a transposition or a product of disjoint transpositions $x_1,x_2\in \Gamma$. Since the operators $\rho(x_i)\in \OR(n)$ commute, the element $\rho(x)$ is an involution and it has eigenvalues $1$ and $-1$. 
The multiplicities $p= \dim \text{Eig}(\rho(x),1)$ and $q= \dim \text{Eig}(\rho(x),-1)$ satisfy the following arithmetic relations  
\begin{align}
\chi(x)=p-q, \quad n=p+q, \quad \det(\rho(x))=(-1)^q, \label{identity1}
\end{align}
from which it can be derived that $2q=n- \chi(x)$. As the group $\Gamma$ is generated by transpositions any any two transpositions are conjugate in $\Gamma$, the representation
$\rho: \Gamma \longrightarrow \OR(n)$ factorizes through $\SO(n)$ precisely when $\det(\rho(x))=1$.
However, identity $2q=n- \chi(x)$ is equivalent to $\chi(x)=n \mod 4$, which shows claim (a).\\[0.3cm] Suppose now that $\rho$ takes values in $\SO(n)$ and let $q=2q'$ for some non-negative integer $q'\in \mathbb{Z}$. The element $\rho(x)\in \SO(n)$ conjugates to the element $s_x\in T_{\SO(n)}$ in the standard maximal torus $T_{\SO(n)}$ given by
$$s_x=\text{diag}(\underbrace{A,...,A}_{q'-times},\ID_p) \in T_{\SO(n)} \ , \quad 
A=\text{exp}_{\so(2)} \left( \begin{array}{cc}
0 & -\pi  \\
\pi & 0  \end{array} \right)=\left( \begin{array}{cc}
-1 & 0  \\
0 & -1  \end{array} \right) \in \SO(2),$$
which can be written in the standard orthonormal basis $(e_1,...e_n)$ of $\mathbb{R}^n$ as $$ \SO(n) \ni s_x=\text{exp}_{\so(n)}(\pi (e_1 \wedge e_2 +...+e_{q'-1}\wedge e_{q'})).$$ Let $\pi:\Spin(n) \longrightarrow \SO(n)$ be the standard covering of $\SO(n)$ and $\cdot: \mathbb{R}^n \times \mathbb{R}^n \longrightarrow \mathbb{R}^n$ be Clifford multiplication. Recall that the differential of the projection $d\pi:\spin(n) \longrightarrow \so(n)$ is a Lie algebra isomorphism that sends $\frac{1}{2}e_i\cdot e_j \in \spin(n)$ to $e_i \wedge e_j \in \so(n)$. The elements in the fiber $\pi ^{-1}(s_x)\subset \Spin(n)$ of $s_x\in T_{\SO(n)}$ are hence of the form 
$$\bar{s}_{x,\varepsilon}=\varepsilon \text{exp}_{\spin(n)}\left( \frac{ \pi}{2} (e_1 \cdot e_{2}+...+e_{q'-1}\cdot e_{q'})\right),$$
where $ \varepsilon\in Z(\Spin(n))=\{1,-1\}$ is a central element in $\Spin(n)$. Direct calculation yields
$$\bar{s}_{x,\varepsilon}^2=\text{exp}_{\spin(n)}^2\left( \frac{\pi}{2} (e_1 \cdot e_{2}+...+e_{q'-1}\cdot e_{q'})\right)=(e_1 \cdot e_{2})^2...(e_{q'-1}\cdot e_{q'})^2=(-1) ^{q'}\ID. $$
Claims (b) and (c) follow, as we have seen that $4q'=n- \chi(x)$.
\epf
\begin{Lem}\label{Pin}
Let $\pi: \Pin(n) \longrightarrow \OR(n)$ be the usual covering and $x\in \SO(n)$ be an element such that $\dim \Eig(x,-1)>0$. Then the following statements hold.
\begin{itemize}
\item[(a)] If $n$ is odd, then the elements in $\pi^{-1}(x)$ are conjugate to each other in $\Spin(n)$.
\item[(b)] The elements in $\pi^{-1}(x)$ are conjugate to to each other in $\Pin(n)$. 
\end{itemize}
 
\end{Lem}
\pf If $n=2m+1$, then the two elements in the fiber of an element $x\in \SO(n)$ with $\dim\Eig(x,-1)>0$ are conjugate to each other in $\Spin(n)$ in view of Lemma 3.9 in \cite{L1}, and \textit{a fortiori} in $\Pin(n)$. So we can suppose that $n=2m$. An element $x\in \SO(n)$ with $\dim\Eig(x,-1)>0$ conjugates to the following element $s_{x}\in T_{\SO(n)}$ in the standard maximal torus $T_{\SO(n)}$ of $\SO(n)$
$$s_{x}=\text{diag}(R(\theta_1), R(\theta_2), ..., R(\theta_m)) \in T_{\SO(n)},$$
where $R(\theta)\in \SO(2)$ is a rotation with angle $\theta\in [0,2 \pi]$ and $\theta_1=\pi$.
The element in the maximal torus $s_{x} \in T_{\SO(n)}$ that corresponds to $x\in \SO(n)$ has fiber $$\pi^{-1}(s_x)=\{\bar{s}_{x,\varepsilon} \in T_{\Spin(n)} \ \ : \ \ \varepsilon \in Z(\Spin(n))\}\subset \Spin(n),$$ where 
$$\bar{s}_{x,\varepsilon}:=\varepsilon\prod_{j=1}^{m}\mu\left(\theta_j\right), \quad  \ \ \mu(\theta_j):= \cos\left(\frac{\theta_j}{2}\right) + (e_{2j-1}\cdot e_{2j} )\sin\left(\frac{\theta_j}{2}\right)\in T_{\Spin(n)} $$
and $\varepsilon\in Z(\Spin(n))$. Conjugation with $e_1 \in \Pin(n)$ exchanges $e_1\cdot e_2\in \Spin(n)$ and $- e_1\cdot e_2\in \Spin(n)$ and fixes $e_j\cdot e_{j+1}\in \Spin(n)$ for any $2 \leq  j<m$, and so it exchanges the elements in the fiber of $s_x\in T_{\SO(n)}$.
\epf
In what follows, we will use the parametrization of conjugacy classes and (complex) irreducible representations of the symmetric group by means of Young diagrams. A standard reference for this material is \cite{sa}.
\bl\label{gen} Two generators of the cyclic group generated by a permutation $z\in \Sym(m)$ are conjugate in $\Sym(m)$.  
\el
\pf We must show that $z^r$ is conjugate to $z$ for any $r\in \mathbb{Z}_{l}^{\times}$, where $l=\ord(z)$. Decompose $z\in \Sym(m)$ as a product of $n$ disjoint cycles $Z_j$ each of length $m_j$ and note that $1=\gcd(r,m_j)$ for any $j\leq n$. Because of the latter condition, the $r$-fold composition of each individual cycle $Z_j$ leads to another cycle of length $m_j$. It follows that $z\in \Sym(m)$ and $z^r\in \Sym(m)$ have the same cycle decomposition and are hence conjugate to each other. 
\epf
\begin{Th}\label{nakayama}
Let $\Gamma=\Sym(m)$ for $6 \neq m \geq 4$ and let $\rho:\Gamma \longrightarrow \OR(n)$ be a faithful irreducible representation of dimension $n=n(\rho)$ whose character $\chi=\chi(\rho)$  fulfills the following conditions
\begin{align}
\chi(x)&= n \mod 8 \label{nakayama2}, \\
\chi(xy) &= n \mod 8 \label{nakayama2.1}, \\
M(z)&:= \frac{1}{l}\sum_{k=0}^{l-1}{(-1)^k\chi(z^k)}>0, \label{nakayama3}
\end{align}
where $x,y \in \Gamma$ are any disjoint transpositions and $z \in \Gamma$ is an arbitrary odd permutation with order $l=l(z)=\ord(z)$. Let $\pi:\Pin(n) \longrightarrow \OR(n)$ be the standard covering homomorphism. The following statements hold. 
\begin{itemize}
 \item[(a)] There is an admissible quadruple $(\Spin(n),\phi_1,\phi_2,\Gamma)$ with $\pi(\phi_i(\Gamma))=\rho(\Gamma)$, provided $n=n(\rho)$ is odd.
 \item[(b)] There is an admissible quadruple $(\Pin(n),\phi_1,\phi_2,\Gamma)$ with $\pi(\phi_i(\Gamma))=\rho(\Gamma)$.
\end{itemize}
\end{Th}
\pf Condition \eqref{nakayama2} implies in view of Lemma \ref{lift} (a) that the representation $(\rho,\Gamma)$ factorizes through $\SO(n)$. Let us consider the canonical central extension $\hat{\Gamma}$ associated with $(\rho,\Gamma)$, i.e. the central extension $\hat{\Gamma}$ making the diagram
\begin{center}
\begin{tikzpicture}
  \matrix (m) [matrix of math nodes,row sep=3em,column sep=4em,minimum width=2em] {
     \hat{\Gamma} &  \Spin(n) \\
     \Gamma & \SO(n) \\};
  \path[-stealth]
    (m-1-1) edge node [right] {} (m-2-1)
            edge node [below] {} (m-1-2)
    (m-2-1.east|-m-2-2) edge node [below] {} node [above] {$\rho$} (m-2-2)
    (m-1-2) edge node [right] {$\pi$} (m-2-2);
\end{tikzpicture}
\end{center}
 commute. Here $\pi:\Pin(n) \longrightarrow \OR(n)$ is the standard covering of $\OR(n)$. \\[0.3cm]
 \textbf{Claim.} Conditions \eqref{nakayama2} and \eqref{nakayama2.1} assure that the central extension $\hat{\Gamma}$ is trivial, i.e. it is isomorphic to $\Gamma \times \mathbb{Z}_2$. 
\pf[Proof of the claim] Recall that isomorphism classes of central extensions of $\Gamma=\Sym(m)$ by $\mathbb{Z}_2$ are in one to one correspondence with elements in the 
second group cohomology group $\text{H}^2(\Gamma,\mathbb{Z}_2)$. It is well-known that the $\mathbb{Z}_2$-(co)homology of the symmetric group stabilizes in degree two. In fact,  $\text{H}^2(\Gamma,\mathbb{Z}_2)$ is isomorphic to the 
Klein Vierergruppe for $m\geq 4$, see the reference to Nakaoka's calculation of the cohomology ring of the symmetric group in \cite{AMI}. The non-trivial elements in $\text{H}^2(\Gamma,\mathbb{Z}_2)$ correspond to group coverings, which we denote by $2 \cdot \Gamma^{\pm}$, and an additional extension 
$L:=2\cdot \Gamma^{+}+2\cdot \Gamma^{-}$. On the other hand $\Gamma=\Sym(m)$ admits a Coxeter presentation \cite{W}
\begin{align}
\Gamma&=\left<t_1,\cdots ,t_{m-1} \ : \ t_i^2=1, \ (t_it_j)^2=1 \ \text{if} \ | i-j|>1, \ (t_it_{i+1})^3=1  \right> \label{csym}.
\end{align} 
The double covers $2 \cdot\Gamma^{\pm}$ also admit a presentation in $m$ generators $z,s_1,...,s_{m-1}$, which collapse to \eqref{csym} once projecting to $\Gamma$, see \cite{sch} page 244.
The relations for $2 \cdot\Gamma^{-}$ read  
\begin{align}
z^2&=1, \label{minus0}\\
s_i^2&=z, \quad  \quad \quad \quad \quad  \quad \    1\leq i\leq m-1 \label{minus1} ,\\
s_{i+1}\cdot s_i \cdot s_{i+1}&= s_{i} \cdot s_{i+1} \cdot s_{i}, \quad \quad  \ \  1 \leq i \leq m-2\label{minus2} ,\\
s_j\cdot s_i&=s_i\cdot s_j\cdot z, \quad \quad  \quad \ \   1\leq i <j\leq m-1, \ \  \ |i-j|\geq2 \label{minus3} ,
\end{align}
whereas for $2 \cdot\Gamma^{+}$
\begin{align}
z^2&=1 ,\\
s_i \cdot z&=z\cdot s_i, \ \ \quad \quad \quad  \quad \  1\leq i\leq m-1,  \\
s_i^2 &= 1, \quad \quad \  \quad \quad  \quad \quad 1\leq i\leq m-1,  \\
s_{i+1}\cdot s_i \cdot s_{i+1}&= s_{i} \cdot s_{i+1} \cdot s_{i}\cdot z, \  \quad 1 \leq i \leq m-2 ,\\
(s_i\cdot s_j)^2&=z, \quad \  \ \ \quad \quad \quad \quad \ 1\leq i <j\leq m-1, \ \  \ |i-j|\geq2 .
\end{align}
Let $x,y\in \Gamma=\Sym(m)$ be disjoint transpositions. Any transposition $x\in \Gamma$ is conjugate to a Coxeter generator $t_i\in \Gamma$, and those generators lift to elements of order two or four depending on the isomorphism type of the canonical extension $\hat{\Gamma}$. 
Let us now consider the order of the lifts of elements $x,xy\in \Gamma$. In case the canonical extension $\hat{\Gamma}$ is trivial, 
these elements lift to involutions. If $\hat{\Gamma}$ is isomorphic to one of the coverings $2\cdot \Gamma^{\pm}$, the order of the lifted elements can be read off from the presentations of such coverings and it is either two or four. Finally, if the canonical extension $\hat{\Gamma}$ is isomorphic to $L$, then one can check that $x\in \Gamma$ lifts to an element of order 4 and $xy \in \Gamma$ lifts to an element of order 2. Table \ref{orders} summarizes these observations. 
\begin{table}[H]
\begin{center}
 \begin{tabular}{ |c || c | c | c| c |} \hline  
Extension & $\Gamma \times \mathbb{Z}_2 $ & $\Gamma^{+} $ & $\Gamma^{-} $ & $L $   \\ \hline
 $\text{Ord}(\bar{x})$ & 2  & 2&4 & 4\\ \hline
 $\text{Ord}(\overline{x y}) $ &2 &4 &4 &2\\ \hline
\end{tabular}
\caption{Order of the lifts of $x\in \Gamma$ and $xy\in \Gamma$.}
\label{orders}
\end{center} 
\end{table} 
Conditions \eqref{nakayama2}-\eqref{nakayama2.1} assure that the elements $x,xy \in \Gamma$ lift to elements of order $2$ in $\hat{\Gamma}$ in accordance with Lemma \ref{lift} (b), and so the extension $\hat{\Gamma}$ must be trivial.
\epf 
In virtue of the last claim, we can choose a homomorphism $\phi_1: \Gamma \longrightarrow \Spin(n)$ lifting the representation $(\Gamma,\rho)$, that is $\pi \circ \phi_1=\rho$, where 
$\pi:\Pin(n) \longrightarrow \OR(n)$ is the standard covering. Put 
$$\phi_2(\gamma)=\eta(\gamma)\phi_1(\gamma), \quad  \gamma \in \Gamma,$$
where $\eta: \Gamma \longrightarrow \mathbb{Z}_2=Z(\Spin(n))$ is the natural non trivial homomorphism.\\[0.3cm]
\textbf{Claim.} The groups $\phi_i(\Gamma)$ are almost conjugate in $\Pin(n)$ for any $n \in \mathbb{N}$, and in $\Spin(n)$ for odd $n\in \mathbb{N}$.
\pf[Proof of the claim] It suffices to check that $\rho(z)$ has $-1$ as an eigenvalue, for any odd permutation $z \in \Gamma $, see Lemma \ref{Pin}.
Consider the cyclic group $\mathbb{Z}_{2l} \subset \Gamma$ generated by an odd permutation $z \in \Gamma$. Denote the restriction of $\rho$ to $\mathbb{Z}_{2l}$ by $\mathcal{Y}|_{\left<z\right>}$ and decompose it into isotypical components 
\begin{align}
V^{\mathcal{Y}|_{\left<z\right>}} &= \bigoplus_{\lambda \in \mathbb{Z}_{2l}}{M_{\lambda} R_{\lambda}} = \bigoplus_{d|2l}{\left(\bigoplus_{  \gcd(\lambda,2l)=d}{M_{\lambda}R_{\lambda}}\right)},\label{sumita}
\end{align} 
where $R_{\lambda}$ is the representation space of the representation of $\mathbb{Z}_{2l}$ with character $\chi_{\lambda}(\cdot)=e^{\frac{2\pi i \lambda (\cdot) }{2l}}$ and $M_{\lambda}$ is its multiplicity in $\mathcal{Y}|_{\left<z\right>}$. Any two generators of $\mathbb{Z}_{2l}=\left<z\right>$ are conjugate in $\Gamma$ by Lemma \ref{gen}. In particular, the action of any element $\sigma \in \Out(\mathbb{Z}_{2l})$ leaves invariant the character of $\mathcal{Y}|_{\left<z\right>}$ and the left-hand side of \eqref{sumita}. In accordance with this, we must have $M_{\lambda}=M_{\sigma(\lambda)}$ for any $\lambda \in \mathbb{Z}_{2l}$. The representation spaces $R_{\lambda}$ in the right-hand side of \eqref{sumita} are in turn permuted according to the action of the element $\sigma \in  \Out(\mathbb{Z}_{2l})$, so we must have $M_{\lambda}=M_d$ for any $\lambda \in \mathbb{Z}_{2l}$ with $\gcd(\lambda,2l)=d$. Consequently, relation \eqref{sumita} becomes 
\begin{align}
V^{\mathcal{Y}|_{\left<z\right>}} &= \bigoplus_{d|2l}{M_d \left(\bigoplus_{  \gcd(\lambda,2l)=d}{R_{\lambda}}\right)}.\quad  \label{sumita2}
\end{align} 
According identity \eqref{sumita2}, the characteristic polynomial of $\rho(z) $ is given by 
$$\det(w-\rho(z))= \prod_{d|2l}\left( \prod_{\gcd(\lambda,2l)=d}{\left(w-e^{\frac{2 \pi i \lambda}{2l}}\right)^{M_d}}\right)=\prod_{d|2l}{\Phi_{\sfrac{2l}{d}}^{M_d}(w)},$$    
where $\Phi_{\sfrac{2l}{d}}$ denotes the $\left(\sfrac{2l}{d}\right)$th cyclotomic polynomial. We conclude that $-1$ is an eigenvalue for $\rho(z)$ precisely when $M_l>0$. 
The explicit formula for the multiplicity $M_l$ is given by \eqref{nakayama3} which is a consequence of the standard orthogonality relations for ireducible characters. 
\epf

\textbf{Claim.} The groups $\phi_i(\Gamma)$ are not conjugate in $\Pin(n)$. Equivalently, there is no automorphism $\sigma\in \Aut(\Gamma)$ and element $g\in \Pin(n)$ such that 
\begin{align}
\Ad_g \phi_1(\sigma(\gamma))& =\phi_2(\gamma), \quad \quad \forall \gamma \in \Gamma. \label{conjug}
\end{align}
 We argue by contradiction. Let us suppose that there is an element $g \in \Pin(n)$ such that \eqref{conjug} holds. Any automorphism of $\Gamma= \Sym(m)$ is inner for $6 \neq m\geq 4$. So, we can choose without loss of generality $\sigma=\ID_{\Gamma}$ and conclude that $\pi(g)\in \OR(n)\subset \U(n)$ centralizes $\rho(\Gamma)$ by projecting relation \eqref{conjug} to $\OR(n)$. Irreducibility of the representation $(\Gamma, \rho)$ tells us that $\pi(g)\in  \U(n)$ is a complex multiple of the identity, and hence $\pi(g)\in \{ \pm \ID \} $. In particular, the element $g \in \left< \omega\right> $, where $\omega = e_1\cdot e_2 \cdot ...\cdot e_n$ is the volume form of the underlying Clifford algebra. Condition \eqref{conjug} tells us then that $\phi_1=\phi_2$, which is absurd.
\epf
The utility of Theorem \ref{nakayama} relies on the fact that its hypotheses can be checked algorithmically: given a representation $\rho$ of $\Sym(m)$ one 
can verify conditions \eqref{nakayama2} and \eqref{nakayama2.1} using the well known Murnaghan-Nakayama rule \cite{sa}, and then create an algorithm to verify condition \eqref{nakayama3} 
manually on all faithful irreducible representations satisfying the former conditions.
\begin{Ex}\label{list}
The table below displays a list consisting of two numbers $n,m \in \mathbb{N}_{\geq 2}$ and a partition $\mathcal{P}=(m_1,m_2,...,m_l)$ of $m$ that defines a representation 
$\rho:\Sym(m) \longrightarrow \OR(n)$ satisfying all hypothesis in Theorem \ref{nakayama}. 

\begin{minipage}{0.2\textwidth}
\begin{flushleft} 
\begin{table}[H]
\begin{center}
 \begin{tabular}{| c || c | c| }  \hline 
$\mathcal{P}$ & $n$ & $m$\\ \hline
(3,2,1) & 16 & 6 \\ \hline
(4,1,1,1,1)& 35& 8\\
(5,2,1)& 64 & \\ \hline
(2, 2, 1, $\cdots$, 1)& 35 & \\
(3, 2, 1, 1, 1, 1, 1)& 160 & \\
(7, 2, 1)& 160 & \\
(5, 4, 1)& 288 & \\
(3, 2, 2, 2, 1)& 288 & \\
(6, 3, 1)& 315 & 10\\
(6, 2, 1, 1)& 350 & \\
(5, 2, 1, 1, 1)& 448 & \\
(3, 3, 2, 1, 1)& 450 & \\
(5, 3, 1, 1)& 567 & \\ \hline
(3, 1, $\cdots$, 1)& 45 & \\
(4, 1, $\cdots$, 1)& 120 & \\
(8, 1, 1, 1)& 120 & 11\\
(7, 4)& 165 & \\
(7, 1, 1, 1, 1)& 210 & \\ \hline
\end{tabular}
\end{center} 
\end{table}
\end{flushleft}
\end{minipage}
\begin{minipage}{1.0\textwidth}
\begin{flushright} 
\begin{table}[H]
\begin{center}
 \begin{tabular}{| c || c | c| }  \hline 
$\mathcal{P}$ & $n$ & $m$\\ \hline
(3, 3, 3, 2)& 462& \\ 
(4, 4, 1, 1, 1)& 825& \\
(4, 2, 2, 1, 1, 1)& 1232 & 11 \\
(6, 3, 1, 1)& 1232& \\ 
(4, 4, 2, 1)& 1320& \\
(4, 3, 2, 2)&1320 & \\ \hline
(4, 1, 1,   $\cdots$ , 1)& 165 &  \\
(3, 2, 1, $\cdots$, 1)& 320&  \\
(9, 2, 1)& 320&  \\
(8, 1, 1, 1, 1)&330 &  \\
(3, 3, 3, 3)& 462&  \\
(6, 1, 1, 1, 1, 1, 1)& 462& 12 \\
(3, 3, 1, $\cdots$, 1)&616 &  \\
(8, 2, 2)&616 &   \\
(6, 5, 1)&1155 &  \\
(5, 5, 2)&1320 &  \\
(3, 3, 2, 2, 2)&1320 & \\
(3, 2, 2, 2, 1, 1, 1)& 1408&  \\ \hline
\end{tabular}
\end{center} 
\end{table}
\end{flushright}
\end{minipage}
\end{Ex}
\begin{remark}
\begin{itemize} 
\item[(i)] The attentive reader might have noticed that condition \eqref{nakayama2} alone  does \textit{not} assure that the 
canonical extension associated with a representation of $\Sym(m)$ is trivial, as claimed in \cite{L2}. For instance, the group $\GL(2,\mathbb{F}_3)$ is an example of a non-trivial extension of $\Sym(4)$ associated with a representation that satisfies condition \eqref{nakayama2} and \textit{not} condition \eqref{nakayama2.1}. The argument in the proof of Proposition 2.6 in \cite{L2} is nonetheless \textit{correct} as the representation of $\Sym(10)$ associated with 
the partition $(2,2,1,\cdots,1)$ fulfills conditions \eqref{nakayama2}-\eqref{nakayama3}. 
\item[(ii)] The representation of $\Gamma=\Sym(6)$ associated with the partition $(3,2,1)$ does not lead to almost conjugate groups of $\Pin(16)$. This is due to the fact that $\Gamma$ is the only symmetric group having non-trivial outer automorphisms, and as there is just one $16$-dimensional representation of $\Gamma$. The pullback of this representation with a non-trivial outer automorphism can be checked to be equivalent to the former one by means of a transformation in $\SO(16)$.
\end{itemize}
\end{remark}

The list in Example \ref{list} suggests that the following conjecture holds. \\[0.3cm]
\textbf{Conjecture.} Let $N(n)$ be the number of equivalence classes of irreducible representations of $\Sym(n)$ satisfying the conditions stated in Theorem \ref{nakayama}. Then, we have $N(n) \geq n$ for any $ n \geq 10.$\\[0.3cm]
Although it is a difficult combinatorial problem to estimate the growth of $N(n)$ against $n\geq 0$, it is not complicated to see that there are infinitely many pairs of almost conjugate subgroups that arise in the spirit of Theorem \ref{nakayama}.
\bp\label{cosito1} 
For each natural number $n \in \mathbb{N}$, there is a pair of almost conjugate subgroups $\Gamma_i(n)$ of $\Spin(m(n))$, where $m(n)=35(2n-1)$, that are not conjugate.
\ep
\pf Let $\Gamma=\Sym(8)$ and let $\rho_1: \Gamma \longrightarrow \SO(V^{35})$ be the 35-dimensional faithful irreducible representation associated with the partition
$(4,1,1,1,1)$. According to Example \ref{list}, conditions \eqref{nakayama2} and \eqref{nakayama3} are met by $\rho_1$. In particular, there is a homomorphism 
$$\phi_1^1:\Gamma \longrightarrow \Spin(V^{35})$$ lifting the representation $\rho_1$. Set 
$\rho_n= \rho \otimes \ID_{\mathbb{R}^n}:\Gamma \longrightarrow \SO(V^{35}\otimes \mathbb{R}^n)$ for an odd number $n\in \mathbb{N}$.
The representation $\rho_n$ satisfies conditions \eqref{nakayama2}-\eqref{nakayama3} and so we can choose a homomorphism 
$$\phi_1^n: \Gamma \longrightarrow \Spin(V^{35}\otimes \mathbb{R}^n)$$ 
lifting the representation $\rho_n$. Define 
$$\phi_2^n(\gamma)= \phi_1^n(\gamma) \eta(\gamma), \quad \gamma \in \Gamma,$$
where $\eta:\Gamma \longrightarrow \mathbb{Z}_2$ is the natural homomorphism, and note that $-1$ is an eigenvalue $\rho_n(\gamma)$ for any odd permutation $\gamma \in \Gamma$. The homomorphisms $\phi_i^n$ are thus almost-conjugate in virtue of Lemma \ref{Pin}. If they where conjugate by some element $y\in \Spin(V^{35}\otimes \mathbb{R}^n)$, then for any $\gamma \in \Gamma $
$$0=[\pi(y),\rho(\gamma) \otimes \ID_{\mathbb{R}^n}]=[A_y \otimes B_y,\rho(\gamma)\otimes \ID_{\mathbb{R}^n}]=[A_y,\rho(\gamma)]\otimes B_y,$$
where $\pi:\Spin(V^{35}\otimes \mathbb{R}^n) \longrightarrow \SO(V^{35}\otimes \mathbb{R}^n)$ is the standard covering and $\pi(y)=A_y\otimes B_y$. Being $\rho$ irreducible, the latter implies that $y= \ID_{\Spin(V^{35})}\otimes \bar{y}$, for some $\bar{y}\in \Spin(\mathbb{R}^{n})$. Consequently $\phi_1^n=\phi_2^n$, which is absurd.
\epf
\section{Isospectral examples}
The manifolds on which we will apply Sunada's criterion are in fact simply connected Lie groups, and as such, they carry a (unique) spin structure. Spin structures of a finite quotient $G_{\Gamma}=
\Gamma \backslash G$ are in one-to-one correspondence with homomorphisms $\varepsilon: \Gamma \longrightarrow \{1,-1\}$. The spinor bundle of $(G_{\Gamma},\varepsilon)$  is in turn given by $ \Sigma_{\varepsilon}( G_{\Gamma})= G \times_{\Gamma} \Spin(n)$, where the action of $\gamma \in \Gamma$
on the first component is the natural one and on $\Spin(n)$ it is by multiplication with the central element $\varepsilon(\gamma)$.  Let $L^2_{\varepsilon}(G_{\Gamma})$ be the space of locally square-integrable complex-valued
$(\Gamma, \varepsilon)$-equivariant functions, i.e. the space consisting of functions $f\in L^2(G)$ such that $f(\gamma g)= \varepsilon(\gamma)f(g)$, for any $\gamma \in \Gamma$ and $  g\in G.$ 
As a Hilbert space $L^2_{\varepsilon}(G_{\Gamma}) \otimes_{\mathbb{C}}  \Sigma_{n}$ can be identified with $L^2(\Sigma_{\varepsilon}(G_{\Gamma }))$ via 
$$ L^2_{\varepsilon}(G_{\Gamma}) \otimes_{\mathbb{C}}  \Sigma_{n} \ni f \otimes s \longmapsto fs \in  L^2(\Sigma_{\varepsilon}(G_{\Gamma })),$$
and Clifford contraction corresponds under this identification to 
$$X \cdot (f \otimes s)= f \otimes (X\cdot s), \quad \quad X \in \gg, \ f\otimes s \in L^2_{\varepsilon}(G_{\Gamma}) \otimes_{\mathbb{C}}  \Sigma_{n} . $$
Ammann and B\"ar observed that if the right $G$-modules $ L^2_{\varepsilon}(G_{\Gamma})$ and $ L^2_{\varepsilon'}(G_{\Gamma'})$ are $G$-isomorphic, then the corresponding Dirac operators are intertwine by the isomorphism of the latter $G$-representation spaces, Theorem 5.1 in \cite{AB}. In paticular, these Dirac operators have the same spectrum. 
If we set the spin structures $(\varepsilon, \varepsilon')\in \Hom(\Gamma,\mathbb{Z}_2) \times \Hom(\Gamma',\mathbb{Z}_2)$ to be the trivial ones, the latter condition is nothing but the assertion that the right regular representations of $G_{\Gamma}$ and $G_{\Gamma'}$ are $G$-isomorphic, which is in turn equivalent to the Sunada condition \cite{Wo}. In particular,  we have the following result.
\bt\label{BA} Let $G$ be a compact simply connected Lie group endowed with a left-invariant metric, and $\Gamma_i \subset G$ be a pair of finite almost conjugate subgroups. Then the finite quotients $G_{\Gamma_1}$ and $G_{\Gamma_2}$ are Dirac isospectral with respect to the  spin structures $\alpha_i(\gamma)=1$ for any $\gamma\in \Gamma_i$ and $i=1,2$.    
  \et 

\subsection{Isospectral nearly K\"ahler pairs}\label{sec2}
Let $G$ be a compact simple Lie group with finite center. The manifold $M=M(G)=G \times G$  admits a left-invariant nearly K\"ahler structure obtained by means of the Ledger-Obata construction \cite{AM}. More explicitly:
identify $M$ with the homogeneous space $L/K$, where $L= G\times G\times G$ and $K=\Delta(G)$ denotes the diagonally embedded $G \hookrightarrow L$. Observe that 
the group $L$ has an obvious 3-symmetry $S \in \Aut (L)$ which stabilizes $K$ and makes $(L,K)$ a 3-symmetric pair. In fact, all homogeneous nearly K\"ahler manifolds are 3-symmetric.\footnote{ We refer to \cite{GM1} for an exhaustive list of homogeneous nearly K\"ahler manifolds.} 
Consequently, 
the Lie algebra $\Delta(\gg) \subset \gl$ admits a complement $\gm \cong T_e M$ invariant under the  
adjoint representation and the infinitesimal action of the 3-symmetry $S_*$. An appropriate rescaled version of the Cartan-Killing form of $L$
\begin{align}
B(X,Y)& =\Tr (\ad(X) \circ \ad(Y)), \quad X,Y \in \gl, \label{NK-metric2}
\end{align}
restricted to the complement $\gm$ defines the nearly K\"ahler metric $g=g_e$ on $M$. The nearly K\"ahler almost complex structure $J$ is in turn an invariant 
tensor characterized on the chosen reductive complement $\gm$ in terms of the $3$-symmetry $S$ as below
\begin{align}
S_{*} &=- \frac{1}{2} \ID + \frac{\sqrt{3}}{2} J\quad \text{on $\gm$}. \label{J}
\end{align}
The triple $(M,g,J)$ defines a homogeneous nearly K\"ahler manifold and we regard 
this structure in the sequel as the \textit{Ledger-Obata space} $M(G)$ associated with $G$. The group of (holomorphic) isometries  of $M=M(G)$ can be calculated using 
that its connected component is $G^3 / \Delta(Z(G))$, where $\Delta(Z(G))$ denotes the diagonally embedded center of $G$, see Theorem 5.3 in \cite{GM}. We include the following lemma for the convenience of the reader and sake of completeness.
\bl\label{isom} Let $M$ be a compact Riemannian manifold and assume that the maximal connected subgroup $G_0$ of $G=\I(M)$ acts transitively on $M$ with stabilizer $K_0$. Then $G = (G_0 \times \Gamma )/L $ as sets, where $\Gamma$ is a subgroup of the quotient of  
$$\Aut(G_0,K_0)=\{ \varphi \in \Aut(G_0) \ : \ \varphi(K_0) \subset K_0\}$$
and the group $\Ad(K_0)$ and $L$ is the kernel of injectivity of the action of $G_0 \times \Gamma$ on $M$. 
\el
\pf Write $M=G_0/K_0$ and let $K$ be the stabilizer of the action of $G$ on $M$. Then, as $G/G_0=K/K_0$, we can write  
$$G= \bigcup_{i}{G_0 \cdot k_i}$$
as sets, where $(k_i)$ is a (finite) system of representatives of $K/K_0$.\\[0.3cm]
On the other hand, the homomorphism $\Ad:K \longrightarrow \Aut(G_0,K_0)$ is injective, see Proposition 1.7 in \cite{sh}.
 We can thus realize the system $(k_i)$ as automorphisms of $G_0$ fixing $K_0$ defined up to elements in $\Ad(K_0)$. 
\epf
The utility of Lemma \ref{isom} relies on the explicit knowledge of the group $\Out(G)$ for a compact semisimple group $G$. For instance, in the complex semisimple case, the latter group corresponds to symmetries of the Dynkin diagram associated with $\gg$.
Direct application of Lemma \ref{isom} to the Ledger-Obata space associated with $G$ leads to the following result.
\bl\label{consequence2} Let $G$ be a compact simply-connected simple Lie group having no non-trivial outer automorphisms and let $M=M(G)$ be the Ledger-Obata space associated with this group. We have
\begin{itemize}
 \item[(a)] The isometry group and the group of holomorphic isometries of $M$ are   
$$\I(M)=G^3/\Delta(Z(G)) \rtimes \Sym(3) \quad \text{and}\quad \I^h(M)=G^3/\Delta(Z(G)) \rtimes \Alt(3),$$
where the action of $\Sym(3)$ on $G^3/\Delta(Z(G))$ is the standard one, and
$\Delta(Z(G))\subset G^3$ denotes the diagonally embedded center.
\item[(b)] Let $\Gamma_i \subset G$ be a pair of non-trivial finite subgroups and identify them with subgroups  
\begin{align}
\overline{\Gamma}_i = \left\{ (\gamma,1,1)\Delta(Z(G))\in \I(M) \ : \ \gamma \in \Gamma_i \ \right\}\subset \I(M). \label{subg}
\end{align}
The groups $\overline{\Gamma}_i$ are conjugate to each other in $\I(M)$ if and only if the groups $\Gamma_i$ are conjugate to each other in $G$.
\end{itemize}
\el
\pf The isometry group and the group of holomorphic isometries of the nearly K\"ahler manifold $(M,g)$ has $L_0=G^3/\Delta(Z(G))$ as maximal connected subgroup. Moreover, the group $L_0$ acts transitively on $M$ with stabilizer $L_{\text{st}}=\Delta(G)/\Delta(Z(G))$.  In virtue of Lemma \ref{isom} we must distinguish the group automorphisms of 
$L_0$ preserving the diagonal $\Delta(G)$ and the metric $g$. Since $G$ has no non-trivial outer automorphisms, the automorphisms of $L_0$ are either inner or permutations of the components. 
The metric $g$ is automatically preserved by those automorphisms, as it is obtained by restriction of the Killing form of $L_0$.  Consequently, the subset $\Gamma\subset \Aut(L_0,L_{\text{st}})$ stated in Lemma \ref{isom} is 
$$\Gamma= \Sym(3)\times N_{G^3}(\Delta(G)),$$ 
where the normalizer $N_{G^3}(\Delta(G))=Z(G)^3\cdot \Delta(G)$ acts on $M$ by conjugation. The kernel of effectivity of the action of $N_{G^3}(\Delta(G))$ on $M$ can be seen to be $ Z(G)^3$.
It follows that, as a set, the group of holomorphic isometries is given by
$$
\I(M)=G^3/\Delta(Z(G)) \times \Sym(3) .$$
Let $\varphi_{[x]} \circ \varphi_{\sigma}\in \I(M)$ be the isometry induced by a coset $([x_1,x_2,x_3],\sigma)$ in the latter quotient. For all $[z_1,z_2,z_3]\in M$, $\sigma,\omega \in \Sym(3)$ and $x,y\in G^3$, we have
\begin{align}
(\varphi_{[x]} \circ \varphi_{\sigma}\circ \varphi_{[y]} \circ \varphi_{\omega})([z])&=
(\varphi_{[(x_1,x_2,x_3)\cdot \sigma(y_1,y_2,y_3)]} \circ \varphi_{\sigma \circ \omega}) ([z]),\nonumber
\end{align}
where $\sigma(x_1,x_2,x_3)=(x_{\sigma(1)},x_{\sigma(2)},x_{\sigma(3)}).$ This shows that the operation in the isometry group $\I(M)$ is given by
$$([x],\sigma)\cdot ([y],\omega)=([x\cdot \sigma(y)],\sigma \circ \omega) \in \I(M).$$
To finish the proof of claim (a) just note that $\Alt(3)\subset \I(M)$ is the permutation group that preserves the almost complex structure \eqref{J}. As for claim (b), let us now
suppose that the groups $\bar{\Gamma}_i \subset \I(M)$ defined in \eqref{subg} are conjugate in $\I(M)$. That is, for any $\gamma_1 \in \Gamma_1$ there is $\gamma_2\in \Gamma_2$ such that 
$$([x],\sigma)\cdot ([\gamma_1,1,1],(1)) = ([\gamma_2,1,1],(1))\cdot ([x],\sigma) $$
for some $([x],\sigma)\in \I(M)$. Evaluation of the latter equality at the trivial coset $\Delta(G)\in M$ yields
$$ \sigma(\gamma_1,1,1)=(\Ad_{x^{-1}}(\gamma_2),1,1) \mod \Delta(G).$$
This is impossible unless $\sigma=(1)\in \Sym(3)$, case in which the groups $\Gamma_i\subset G$ are conjugate to each other in $G$.
 \epf
We are now in position to give a proof for the main result of this paper.
\bt  \label{isospectral} There is a strictly increasing sequence of numbers $(d_n)_{_{n\in \mathbb{N}}}$ such that for each $n\in \mathbb{N}$ there is a pair $M^{d_n}_1$ and $M^{d_n}_2$  of non-isometric nearly K\"ahler manifolds that are isospectral for the Dirac and the Hodge-Laplace operator $\Delta^k$ for $k=0,1,...,\dim(M)$.
\et
\pf For each $n\in \mathbb{N}$ choose an odd number $d_n$ so that $G=\Spin(d_n)$ satisfies the hypotheses in Theorem \ref{nakayama} or Proposition \ref{cosito1}. In particular, $G$ admits almost conjugate subgroups $\Gamma_i$ that are not conjugate to each other. Denote the Ledger-Obata space associated with $G$ by $M=M(G)$. The chosen pair of almost conjugate groups $\Gamma_i$ of $G$ act freely on $M(G)=G \times G$ by left multiplication in the first component and yield an isospectral pair $M_i$. Lemma \ref{consequence2} implies that the quotients $M_i$ are not isometric and hence the claim for the Laplace operator. The claim for the Dirac operator follows from Theorem \ref{BA} endowing the quotients $M_i$ with the spin structures $\alpha_i(\gamma)=1$ for any $\gamma \in \Gamma_i$, where $i=1,2$. 
\epf
\begin{remark} Almost conjugate finite groups $\Gamma_i$ constructed by means of Theorem \ref{nakayama} yield 
  explicit examples of isospectral good Riemannian orbifolds with different \textit{maximal isotropy orders}, see Corollary 2.6 in \cite{RSW} and good isospectral spin orbifolds $\Gamma_i \backslash S^{2n}$. 
\end{remark}
\subsection{Dimension six}
Sunada isospectral quotients are generic in the compact setting, see \cite{P}, and so it is natural to ask for the existence of Sunada pairs of nearly K\"ahler manifolds in dimension six. The aim of this section is to show the following result
\bp\label{NOT} If $M_{\Gamma}$ and $M'_{\Gamma'}$ are non-isometric Laplace isospectral locally homogeneous nearly K\"ahler manifolds in dimension six with $\Gamma \subset \I_0^h(M)$ and $\Gamma' \subset \I_0^h(M')$,  then  $M'$ and $M$ are holomorphically isometric to $S^3\times S^3$ endowed with its 3-symmetric nearly K\"ahler structure. Furthermore,
	 Sunada isospectral pairs $M_{\Gamma}$ and $M_{\Gamma'}$ with $M=S^3\times S^3$ and $\Gamma, \Gamma' \subset \SU(2) \times \SU(2) \times \{\ID\}$ are holomorphically isometric.
	\ep
 In order to show this statement, we will make use of the following formula to compute the volume of a compact Lie group with respect to the Haar measure induced by an invariant metric: 
  \bp\label{volume} Let $G$ be a compact, connected Lie group, let $T \subset G$ be a maximal torus and $W(G)$ be the Weyl group of $G$.
  The volume of $G$ with respect to a metric induced by an Ad-invariant scalar product $B$ on $\gg$ is given by 
  \begin{align}
   \vol(G,B)= 2^{\dim(T)} | W(G)|\pi^{\frac{\dim(T) +\dim(G)}{2}} \left( \sqrt{|\det(b^{-1}(\e_{\mu},\e_{\nu}))|} \  e^{-\frac{1}{4} \Delta}|_0 \delta_{\gg}\right)^{-1}, \label{vol}
  \end{align}
  where $ (\e_{\mu})$ is a basis of the dual weight lattice $P(G)^*\subset \gt^*$, and $\delta_{\gg}$ denotes the product of the product of the roots of $\gg ^{ \mathbb{C}}$ and $b^{-1}$ is the scalar product of $\gt^{*}$ induced by $B$. $\Delta$ is the unique second order linear differential operator with constant coefficients acting on the polynomial algebra $\mathbb{C}[\varepsilon_1,...,\varepsilon_m]^{\Sym(m)}$ in the simple roots $ \varepsilon_1,...,\varepsilon_m \in \gt^{*}$ of $\gg^{\mathbb{C}}$, such that $-\Delta\alpha^{2}=2b^{-1}(\alpha,\alpha)$ and $\Delta\alpha^{2k+1}=0$ for any functional $\alpha=\varepsilon_1+...+\varepsilon_j \in  \gt^*$ and $j\leq m$.  
\ep
\pf The formula is clear for a torus with respect to the metric induced by $B$. That is, 
$$\vol(T,B)=(2\pi)^{\dim(T)}\ |\det(b^{-1}(\e_{\mu},\e_{\nu}))|^{-1/2}.$$
Define a function $f$ on $\gg$ by 
$ f(X)=e^{-B(X,X)}$
and observe that 
\begin{align}
 \pi^{\frac{\dim(G)}{2}}&= \int_{\gg} f(X)d_B(X)= \int_{\gg^{reg}}f(X) d_B(X), \label{f0}
\end{align}
where $\gg^{reg}\subset \gg$ is the generic set consisting of regular elements, and $d_B$ denotes the Riemannian volume of the metric $B$.
Choose an open Weyl chamber $\mathcal{C}$ and consider the map 
\begin{align} 
\gg^{reg} \ni X \longmapsto y(X) \in \mathcal{C}, \label{prob}
\end{align}
where $y(X)\in \Ad_{G}(X)\cap \mathcal{C}$. The metric $B$ on $\gg$ induces metrics on $\mathcal{C}$, $\gt$, $\gg^{reg}$ and $\gt^{reg}$, and the map defined in \eqref{prob} becomes a Riemannian submersion with respect to these metrics. 
 Let us denote by $d_b$ the Riemannian volume of the metric $b$ on $\gt$ and by $d_{E_X}$ the one on the adjoint orbit of a regular element $X\in \gg^{reg}$ relative to the submersion \eqref{prob}. Integration over the submersion \eqref{prob} yields 
\begin{align}
\int_{\gg^{reg}}f(X) d_B(X)& =\frac{1}{|W(G)|} \int_{\gt^{reg}} \left(\int_{\Ad_{G}(X)} f(Y) d_{E_X}(Y) \right)d_b(X) \nonumber \\
&=  \frac{1}{|W(G)|}\int_{\gt} f(X) \left(\int_{\Ad_{G/T}(X)}  d_{E_X}(Y) \right)d_b(X) , \label{f1} 
\end{align}
where in the last equality we have used the Ad-invariance of $f(\cdot)$ and the fact that $T$ is a maximal torus whose Lie algebra contains $X \in \gt^{reg}$. On the other hand, the orbit map $$\Ad_{G/T}:gT \longmapsto \Ad_g(X)$$ provides an isometric immersion of $(G/T, \Ad_{G/T}^*B )$ into $(\gg^{reg},B)$. Comparison of the volumes of the flag manifold $G/T$ relative to its submersion metric $B$ and $\Ad_{G/T}^*B$ gives
$$\vol(G/T, \Ad_{G/T}^*B )=\det(\ad_X|_{\gt^{\perp}})\vol(G/T,B)=\delta_{\gg}(X)\frac{\vol(G,B)}{\vol(T,B)},$$
where $\delta_{\gg}$ is the product of the roots of $\gg ^{\mathbb{C}}$. As a consequence,
we have  
\begin{align}
\int_{\Ad_{G/T}(X)}d_{E_X}(Y)=  \frac{\vol(G,B)}{\vol(T,B)}\delta_{\gg}(X), \label{ff}
\end{align}
and identities \eqref{f0},\eqref{f1} and \eqref{ff} give
\begin{align}
 \pi^{\frac{\dim(G)}{2}}&= \int_{\gg} f(X)d_B(X)= \frac{\vol(G,B)}{\vol(T,B)} \ \frac{1}{|W(G)|} \int_{\gt}{f(X) \delta_{\gg}(X) d_{b}(X) }. \label{f3}
\end{align}
Formula \eqref{vol} is then a consequence of identity \eqref{f3} and the following claim.\\[0.3cm]
\textbf{Claim.} The following formula holds for any polynomial $P \in \mathbb{C}[\varepsilon_1,...,\varepsilon_m]^{\Sym(m)}$ in the simple roots $\varepsilon_1,...,\varepsilon_m \in \gt^{*}$.     
\begin{align}
\frac{1}{\sqrt{\pi^{\dim(T)}}}\int_{\gt}{f(X) P(X) d_b(X) }&= e^{-\frac{1}{4}\Delta}|_0 P. \label{form}
\end{align}
The standard polarization argument \cite{gt}, tells us that it suffices to verify formula \eqref{form} for $P=\alpha^k$ for any $k\in \mathbb{N}$ fixed and any functional $ \alpha= \varepsilon_1 + ... +\varepsilon_j \in \gt^*$, where $\varepsilon_i \in \gt^*$ is a simple root of $\gg^{\mathbb{C}}$ and $j \leq m$. First observe that 
\begin{align}
\sum_{k\geq 0}{ \frac{t^k}{k!}} & \cdot \frac{1}{\sqrt{\pi^{\dim(T)}}} \int_{\gt}{e^{-b(X,X)}\alpha^k(X) d_b(X)}=  \frac{1}{\sqrt{\pi^{\dim(T)}}} \int_{\gt}{e^{-b(X,X)}\cdot e^{t \alpha(X)} d_b(X)} \nonumber\\
&=  e^{\frac{t^2}{4}b^{-1}(\alpha,\alpha)}= \sum_{k\geq 0 }{\frac{t^{2k}}{k!} \cdot \frac{b^{-1}(\alpha,\alpha)^{k}}{4^k} }.\nonumber
\end{align}
The left hand side of identity \eqref{form} vanishes for $P=\alpha^k$ with odd $k\in \mathbb{N}$ and 
$$\frac{1}{\sqrt{\pi^{\dim(T)}}} \int_{\gt}{e^{-b(X,X)}\alpha^{2k}(X) d_b(X)}=  \frac{b^{-1}(\alpha,\alpha)^{k}}{4^k}\cdot  \frac{(2k)!}{k!}.$$
The right hand side of \eqref{form} vanishes for $P=\alpha^k$ with odd $k\in \mathbb{N}$ according to the definition of the operator $\Delta$. Moreover, we have
\begin{align}
e^{-\frac{1}{4}\Delta}|_0 \alpha^{2k}&= \sum_{k\geq l \geq 0}{\frac{1}{4^l} \frac{b^{-1}(\alpha,\alpha)^l \alpha^{2k-2l}(0)}{l!} \cdot \frac{(2k)!}{(2k-2l)!} } =\frac{b^{-1}(\alpha,\alpha)^{k}}{4^k}\cdot  \frac{(2k)!}{k!}. \nonumber 
\end{align}
The claim follows.
\epf
\begin{remark}
The idea of using Weyl integration's formula to reduce the problem of computing the volume of a compact Lie group $G$ to the computation of the volume of the associated flag variety appears in Macdonald's work \cite{mac}. Indeed in part (2) of his paper, Macdonald obtained a fomula that is essentially identity \eqref{f3}. The rest of the proof of Proposition \ref{volume} can be thought as a simplification of his argument to compute the integral given in the left hand side of \eqref{form}.
\end{remark}
	\bp\label{VOL} The volume of the four homogeneous six dimensional nearly K\"ahler manifolds with respect to the unique invariant metric $g\in \Sym^2(TM)$ with scalar curvature $\kappa$ is given by Table \ref{volume1}.
\begin{table}[H]
\begin{center}
 \begin{tabular}{ |c || c | c | c| c |} \hline  
 $(M^6,g)$ & $S^6$ & $\mathbb{C}P^3$ & $S^3 \times S^3$ & $F(1,2)$\\ \hline
$\vol(M^6)$& $\left(\frac{30}{\kappa}\right)^3 \frac{16\pi^3 }{15}$&$\left(\frac{30}{\kappa}\right)^3 \frac{\pi^3 }{6}$ & $\left(\frac{30}{\kappa}\right)^3 \frac{8\pi^4 }{81\sqrt{3}}$ 
& $\left(\frac{30}{\kappa}\right)^3 \frac{\pi^3 }{2}$ \\ \hline
\end{tabular}
\caption{Volumes of homogeneous $(M^6,g)$ with $\text{scal}(g)=\kappa$.}
\label{volume1}
\end{center} 
\end{table}
\ep\pf The nearly K\"ahler metric $g  \in \Sym^2(TM)$ on $M=G/K$ with $G=\I_0^h(M)$, is induced by an invariant scalar product $B$ on $\gg$. Furthermore, the fibers of the Riemannian submersion $\pi:G \longrightarrow G/K$ have the same volume with respect to the left invariant metric induced by $B$ on $G$, and so 
$$\vol(M,g)=\vol(G, B)\vol(K, B|_{\gk \times \gk})^{-1}.$$
 Consequently, we can use Proposition \ref{volume} to calculate the volume of $M$ with respect to a chosen normal metric. We choose our metric such that the scalar curvature of $M$ is $30$, Lemma 5.4 in \cite{MS}, and perform the explicit calculations of the volumes just for $\mathbb{C}P^3$, $F(1,2)$ and $S^3 \times S^3$. 
\begin{itemize}
\item[(i)] Let us first consider the flag manifold $F(1,2)=\SU(3)/T$ with the normal metric induced by the invariant scalar product 
$$B:\su(3) \times \su(3) \longrightarrow \mathbb{R}, \quad \quad  B(X,Y)=-\frac{1}{2}\tr_{\mathbb{C}^3}(XY),\quad X,Y \in \su(3).$$
Fix a maximal torus of $\SU(3)$ by declaring its Lie algebra to be
$$\gt=\{\diag(iX_1,i X_2,iX_3) \ : \ X_i \in \mathbb{R}, \ X_1 +X_2 +X_3=0\}=\text{span}(H_i \ : \ i\leq 3),$$
where  $H_1=\diag(0,i,-i)$, $H_2=\diag(-i,0,i)$ and $H_3=-H_1-H_2=\diag(i,-i,0)$. Eventually,  in the base $(H_1,H_2)$ and $(\varepsilon_1,\varepsilon_2)$, we get
$$b(H_{\mu},H_{\nu})=\frac{3}{2}\left(\delta_{\mu\nu}-\frac{1}{3}\right) \quad \text{and}\quad \ b^{-1}(\varepsilon_{\mu},\varepsilon_{\nu})=2\left(\delta_{\mu\nu}-\frac{1}{3}\right),$$
where $\varepsilon_j \in \gt^*$ are the functionals $\varepsilon_j(\diag(iX_1,iX_2,iX_3))=X_j$ 
and $ b^{-1}:\gt^* \times \gt^* \longrightarrow \mathbb{R}$ is the scalar product dual to $b=B|_{\gt \times \gt}$.
The roots of $ \su(3)$ are 
$$\pm i(\varepsilon_1-\varepsilon_2),\pm i(\varepsilon_3-\varepsilon_1), \pm i(\varepsilon_2-\varepsilon_3),$$ and a straightforward calculation reveals that
$$-\Delta= \frac{4}{3}\left( \frac{\partial^2}{\partial^2\varepsilon_1^2}+\frac{\partial^2}{\partial^2\varepsilon_2^2}\right)
-\frac{4}{3}\frac{\partial^2}{\partial \varepsilon_1 \partial \varepsilon_2},$$
where $\Delta$ is the operator introduced in Proposition \ref{volume}.
The polynomial $\delta_{\su(3)}$ is
$$\delta_{\su(3)}=(\varepsilon_1-\varepsilon_2)^2(2\varepsilon_2+\varepsilon_1)^2(2\varepsilon_1+\varepsilon_2)^2 =4 \varepsilon_1^6 + 12 \varepsilon_1^5\varepsilon_2-3\varepsilon_1^4 \varepsilon_2^2 
-26\varepsilon_1^3 \varepsilon_2^3-3\varepsilon_1^2\varepsilon_2^4+12\varepsilon_1\varepsilon_2^5 +4\varepsilon_2^6,$$
and hence 
$$ e^{-\frac{1}{4} \Delta}|_0 \delta_{\su(3)}= \frac{1}{6\cdot 3^3}\left( \frac{\partial^2}{\partial \varepsilon_1^2}+\frac{\partial^2}{\partial \varepsilon_2^2}-\frac{\partial^2}{\partial \varepsilon_1\partial \varepsilon_2}\right)^3 \delta_{\su(3)}(0) = \frac{1944}{162}= 12.$$
Using the fact that $|W(\SU(3))|=|\Sym(3)|=6$, formulae  \eqref{vol} and \eqref{f3} yield 
$$ \vol(F(1,2),B)=|W(\SU(3))|(e^{-\frac{1}{4} \Delta}|_0 \delta_{\su(3)})^{-1}\pi^{\frac{\dim(\su(3))-\dim(\gt)}{2}}=\frac{1}{2}\pi^3.$$
\item[(ii)] 
Let us now consider $\mathbb{C}P^3=\Sp(2)/\U(1)\times \Sp(1)$ with the submersion metric induced by 
$$B(X,Y)=-\frac{1}{2}\Real \tr_{\mathbb{H}^2}\left(XY\right) \quad \quad X,Y\in \mathfrak{sp}(2).$$
Fix a maximal torus $T$ in $\Sp(2)$ by declaring 
$\gt=\left\{ \text{diag}(i\theta_1,i\theta_2)\ : \ \theta_i \in \mathbb{R}\right\}$ 
and observe that 
$$b=B|_{\gt\times \gt}=\frac{1}{2}(\theta_1 \otimes \theta_1+ \theta_2 \otimes \theta_2),$$
where $\theta_i \in \gt^{*}$ are the functionals 
$$\theta_i(\text{diag}(i\tilde{\theta}_1,i\tilde{\theta}_2))=\tilde{\theta}_i, \quad \text{diag}(i\tilde{\theta}_1,i\tilde{\theta}_2) \in \gt.$$
Consequently, we have 
\begin{align}
b^{-1}(\theta_1,\theta_1)=b^{-1}(\theta_2,\theta_2)&=2 \quad \text{and} \quad b^{-1}(\theta_1,\theta_2)=0.\label{co}
\end{align}
On the other hand, the roots of $\gsp(2)$ are $\pm 2i \theta_1, \pm 2i\theta_2, \pm i\theta_1 \pm i \theta_2$ and 
$$\delta_{\gsp(2)}=16\ \theta_1^2\theta_2^2(\theta_1 ^2 -\theta_2^2)^2.$$ From the identities in \eqref{co}, we see that the operator described in Proposition \ref{volume} is in turn given by 
$$\Delta=4 \left( \frac{\partial^2}{\partial \theta_1}+\frac{\partial^2}{\partial \theta_2}\right).$$
After a short calculation, we get
$e^{\frac{-1}{4}\Delta}\delta_{\gsp(2)}|_0=192.$
In view of this and the fact that $|W(\Sp(2))|=8$, formula \eqref{vol} gives
$\vol(\Sp(2),B)=\frac{\pi^6}{12}.$
We also note that $$ \vol(\U(1)\times \Sp(1),B|_{\gu(1)\oplus \gsp(1)})=\frac{\pi^3}{2},$$ and so
$ \vol(\mathbb{C}P^3,B)=\frac{\pi^3}{6}.$
\item[(iii)] Let us now consider $S^3\times S^3=\SU(2)^3/\Delta(\SU(2))$, where $\Delta(\SU(2))$ denotes the diagonally embedded $\SU(2)$ in $\SU(2)\times \SU(2) \times \SU(2)$, endowed with the metric 
$$B(X,Y)=-\frac{1}{3}\tr(XY), \quad \quad X,Y\in \su(2) \oplus \su(2) \oplus \su(2).$$ 
We have $B(\Delta H,\Delta H)=2$, and $B^{-1}(\Delta x, \Delta x)=\frac{1}{2}$, where 
$$\Delta H=(H,H,H)\in \Delta(\su(2)) \quad , \quad H=\left( \begin{array}{cc}
i & 0  \\
0 & -i  \end{array} \right) \in \su(2)$$ 
and $\Delta x\in (\Delta\su(2))^*$ is the unique element such that $\Delta x(\Delta H)=1$. The polynomial $\delta_{\Delta(\su(2))}=4 \Delta^2 x$, and so the volume of $\Delta(\SU(2))$ is 
\begin{align}
\vol(\Delta(\SU(2)), B|_{\Delta(\su(2))\times \Delta(\su(2))})&=32\sqrt{2}\pi ^2.
\end{align} 
Similarly, one calculates that
$$\vol(\SU(2) \times \SU(2)\times \SU(2), B)^{1/3}= \frac{8\sqrt{2}\pi^2}{3\sqrt{3}},$$
and so 
$$\vol(S^3 \times S^3,B)=\frac{32\pi^4}{81\sqrt{3}}.$$
\end{itemize}
Table \ref{volume1} summarizes the results.
\epf
 \pf[Proof of Proposition \ref{NOT}]  We can distinguish the Laplace spectrum of two locally homogeneous six dimensional nearly K\"ahler manifolds $(M_{\Gamma},g)$ and $(M'_{\Gamma'},g')$ by their heat kernel coefficients, i.e. the coefficients of the asymptotic expansion of their Laplacians $\Delta^0_g$ and $\Delta^0_{g'}$ acting on functions. In fact, heat invariants up to order 1 of $(M_{\Gamma},g)$ and $(M'_{\Gamma'},g')$ are given by, see Theorem 4.8.18 in \cite{Gi}:
$$a_0^0(g')=\vol(M'_{\Gamma'}), \quad a_1^0(g')=\frac{1}{6}\int_{M'_{\Gamma'}}{ \kappa(M'_{\Gamma'}) d\vol(M'_{\Gamma'})}$$
and
$$a_0^0(g)=\vol(M_{\Gamma}), \quad a_1^0(g)=\frac{1}{6}\int_{M_{\Gamma}}{ \kappa(M_{\Gamma}) d\vol(M_{\Gamma})},$$
where $\kappa(M'_{\Gamma'})$ and $\kappa(M_{\Gamma})$ denote the scalar curvature of $(M'_{\Gamma'},g')$ and $(M_{\Gamma},g)$ respectively. Since nearly K\"ahler six manifolds are Einstein, isospectrality of $(M_{\Gamma},g)$ and $(M'_{\Gamma'},g')$ implies $\kappa(M)=\kappa(M_{\Gamma})=\kappa(M'_{\Gamma'})=\kappa(M')$ and 
 $$\frac{\vol(M)}{\vol(M')}=\frac{|\Gamma|}{|\Gamma'|} \in \mathbb{Q}. $$
The first claim is a consequence of this reasoning together with Proposition \ref{VOL} and Proposition 1.1 in \cite{CV}, whereas the second follows from Proposition \ref{SU2}.  
\epf

 
\end{document}